\documentclass[12pt]{article}
\usepackage{amsmath}
\usepackage{amssymb}
\usepackage{enumitem}
\usepackage{amsthm}
\usepackage{xcolor}
\usepackage{cancel}
\usepackage[left=2.4cm,right=2.4cm,top=2.4cm,bottom=2.4cm]{geometry}
\usepackage[colorinlistoftodos]{todonotes}
\usepackage[backref]{hyperref}

\newtheorem{theorem}{Theorem}
\newtheorem{lem}[theorem]{Lemma}

\newtheorem{lemma}[theorem]{Lemma}

\def \btm {\begin{theorem}}
\def \etm {\end{theorem}}

\newtheorem{prop}[theorem]{Proposition}
\newtheorem{rmk}[theorem]{Remark}
\newtheorem{cor}[theorem]{Corollary}

\newtheorem{constr}[theorem]{Construction}

\def \bpf {\begin{proof}}
\def \epf {\end{proof}}

\def \nin {\noindent}
\def \bsk {\bigskip}

\def \smin {\setminus}

\def \eee {\mathbb{E}}
\def \ppp {\mathbb{P}}
\def \cB {\mathcal{B}}

\def \blm {\begin{lem}}
\def \elm {\end{lem}}

\def \cH {\mathcal{H}}

\def \dist {\mathrm{dist}}

\def \ex {{\rm ex}}

\def \ind {{\rm \mbox{\rm -}ind}}    
\newcommand{\turind}[2]{\ex(n; \,#1, \,#2\ind)}

\begin{document}

\title{Graphs with constant links\\ and induced Turán numbers
}
\author{Yair Caro\thanks{~Department of Mathematics, University of Haifa-Oranim, Tivon 36006, Israel} \and Adriana Hansberg\thanks{~Mathematics Institute, UNAM Juriquilla, Quer\'etaro, Mexico} \and Zsolt Tuza\thanks{~HUN-REN Alfr\'ed R\'enyi Institute of Mathematics, Budapest, Hungary} $^,$\thanks{~Department of Computer Science and Systems Technology, University of Pannonia, Veszpr\'em, Hungary}}
%\date{\small Latest update on \version}
\date{\small Latest \LaTeX\ run on \today}
\maketitle

\begin{abstract}
A graph $G$ of constant link $L$ is a graph in which the neighborhood of any vertex induces a graph isomorphic to $L$. Given two different graphs, $H$ and $G$, the induced Turán number $\turind{H}{G}$ is defined as the maximum number of edges in an $n$-vertex graph having no subgraph isomorphic to $H$ and no copy from $G$ as an induced subgraph. %Given two disjoint families of graphs, $\mathcal{P}$ and $\mathcal{Q}$, the induced Turán number $\turind{\mathcal{P}}{\mathcal{Q}}$ is defined as the maximum number of edges in an $n$-vertex graph having no member from $\mathcal{P}$ as a subgraph and no member from $\mathcal{Q}$ as an induced subgraph.

Our main motivation in this paper is to establish a bridge between graphs with constant link
and induced Tur\'an numbers via the class of $t$-regular, $k$-uniform (linear) hypergraphs of girth at least $4$, as well as to present several methods of constructing connected graphs with constant link.  

We show that, for integers $t \geq 3$ and $k \geq 3$, $\turind{C_k}{K_{1,t}} \leq (k - 2)(t - 1)n/2$ and that equality holds for infinitely many values of $n$. This result is built upon the existence of graphs with constant link $tL$ with restricted cycle length, which we prove in another theorem. More precisely, we show that, given a graph $F$ with constant link $L$ and circumference $c$, then, for all integers $t \geq 2$ and $g > c$, there exists a graph with constant link $tL$ which is free of cycles of length $l$, for all $c < l < g$. We provide two proofs of this result using distinct approaches. We further present constructions of graphs with constant links $tL$, $t \geq 2$, and restricted cycle length based on Steiner Systems. Finally, starting from a connected graph of constant link $tL$, for $t \geq 2$, having order $n$ and restricted cycle lengths, we provide a method to construct an infinite collection of connected graphs of constant link $tL$ that preserves the cycle length restriction, and whose orders form an arithmetic progression $qn$, $q \geq 1$.   
\bsk

\nin
\bsk

\nin
\end{abstract}

\section{Introduction}

Our main motivation in this paper is to establish a bridge between graphs with constant link and induced Turán numbers via the class of regular, uniform (linear) hypergraphs of girth at least $4$, as well as to contribute with several methods of constructing connected graphs with constant link. In order to explain this more precisely, we will present a short overview on both, link graphs and induced Turán numbers.

\subsection{Link graphs}

Let $G$ be a finite undirected graph. The \emph{link} of a vertex $v$ is the subgraph $L(v):=G[N(v)]$ induced by the open neighborhood $N(v)$ of $v$. If the links $L(v)$ are isomorphic to the same graph $H$ for all $v$, then $G$ is said to have \emph{constant link} $H$---or, $G$ is a \emph{constant-link graph}, if the actual structure of $H$ is unimportant or understood---and $H$ is called a \emph{link graph}. Zykov asked in 1963 \cite{Zyk63}, what is sometimes called the Trahtenbrot--Zykov problem, whether, given a finite graph $H$, there exists a graph $G$ (not necessarily finite) with constant link $H$. In general, as an algorithmic problem, this question is undecidable \cite{Bul73}. It is unknown whether the problem restricted to finite $G$ is also undecidable.

More generally, one wants to classify all graphs that are link graphs. Early works discussed necessary conditions for a graph to be a link graph \cite{BroCon73,BroCon75} or a link tree \cite{BlaHM80}; and all link graphs on at most six vertices were identified \cite{Hal85}. Since then, over 100 papers appeared discovering many new link graphs. Most of the methods used for the construction of link graphs are based on deep structures of corresponding groups and incidence structures coming from combinatorial geometry. As we have no intention here to survey the literature on link graphs, we refer the reader to the papers \cite{BlaHM80, Brou23, Brou95, BroCon73, BroCon75, LaMu13, Toma89, Vog86, Weet94-1, Weet94-2} where many further details and more references can be reached. 
See also \cite{LaPiVi19, YuWu21, Zel86} for related problems.

We mention here that, except for sporadic cases, there are only few papers concerning link graphs of the form $tH$, i.e.\ $t$ vertex-disjoint copies of $H$, or, more generally, link graphs which are unions of smaller graphs. The most studied and mentioned cases are $tK_2$, linear forests, and graphs of the form $H \cup tK_1$ where $H$ is a link graph \cite{Fro89}. An explicit  construction for the case that $H_1 \cup H_2$ is a link graph provided $H_1$ and $H_2$ are link graphs is given in Theorem~U of \cite{Hal85}, aimed to show that there are infinitely many non-isomorphic graphs with $H_1 \cup H_2$ as a link graph. However, this construction is otherwise unnecessarily involved, as the reader can see by the following Cartesian product construction (which seems to have been observed before several times) by which it is easily seen that, for $1 \leq j \leq t$, if $G_j$ is a graph with constant  link $H_j$, then the Cartesian  product $G = G_1\Box G_2\Box \cdots \Box G_t$ is a graph with constant link $H = \cup_{j = 1}^tH_j$. On the other hand, the Cartesian product of a graph $G$ may have many cycles not present in $G$; for example, $K_p\Box K_p$ is pancyclic for any $p \geq 3$. The same happens in the construction suggested in \cite{Hal85} as no care was given to the appearance of short cycles. We are unaware of any constructions of link graphs with attention given to the appearance of short cycles. This point is crucial in establishing the connection between link graphs and induced Turán numbers.

Let us note also that there are link graphs $H$ that appear as link graphs of only finitely many connected graphs.
A trivial example is $K_k$.
More involved, the Petersen graph is known to be the link graph of only three connected graphs \cite{Hal80}. Hence, in general, it makes sense to consider the construction of infinitely many connected graphs with link $tH$, $t \geq 2$. Another question in the opposite direction with this regard is the existence of a graph $H$ such that $tH$ is a constant link for some $t \geq 2$ but $jH$ is not a constant link for $1 \leq j < t$.

\subsection{Induced Turán numbers}

Let $\mathcal{P}$ and $\mathcal{Q}$ be two disjoint families of graphs and define $\turind{\mathcal{P}}{\mathcal{Q}}$ as the maximum number of edges in an $n$-vertex graph having no member from $\mathcal{P}$ as a subgraph and no member from $\mathcal{Q}$ as an induced subgraph. Induced Turán problems were explicitly introduced by Loh, Tait, Timmons and Zhou \cite{LTTZ18}, with the main result stating that, for any graph $H$ and integers $t \geq s \geq 1$, $\turind{H}{K_{s,t}} = O (n^{2-\frac{1}{s}})$. Since then, the subject has received much interest. We refer the reader to the following papers for further information: \cite{EGM19} concerning sharp asymptotic bounds for $\turind{C_{2k+1}}{K_{s,t}}$, \cite{NTT18} concerning $\turind{H}{K_{s,t}}$ (for any graph $H$) using spectral methods and improving upon \cite{LTTZ18}, and \cite{Ill21-1, Ill21-2} with respect to the structure of graphs without induced $K_{2,t}$, the asymptotics of $\turind{H}{F}$ when $F$ and $H$ have chromatic number at least $3$, and asymptotically sharp bounds on $\turind{H}{F}$ for any two graphs $H, F$, with exception of $H$ being bipartite or $F$ an empty graph or a complete bipartite graph.

Asymptotic bounds were established in many cases when $t \geq s \geq 2$, and, in particular, for $\turind{C_k}{K_{2,t}}$. However, in the case that $s = 1$, namely $\turind{H}{K_{1,t}}$, a linear upper bound is easy to establish, but the exact value was left unsolved even for cycles. In a companion paper (under preparation), we consider $\turind{H}{K_{1,t}}$ in more generality. However, in the specific case of $\turind{C_k}{K_{1,t}}$, an exact value can be obtained using a connection to link graphs, thus also answering a request from \cite{LTTZ18} to get better estimate on $\turind{\mathcal{Q}}{\mathcal{P}}$, which we achieve in the corresponding cases.

\subsection{Regular, uniform hypergraphs of large girth}

Problems concerning the existence of $t$-regular $k$-uniform hypergraphs with large girth were the subject of much research since the initial work done by Erd\H{o}s and Sachs concerning the case of graphs, that is when $k =2$ \cite{ErSa63}. The main problems were of two types. The first was the existence of such hypergraphs and the second was their lowest possible order (a problem inspired by the notion of the well known Moore bound). For the purpose of this paper, the following results suffice.

\btm[\cite{FLSUW95}] \label{t:existence_Furedi_etal}
For every $k, t, g \geq 2$, there is a $k$-uniform $t$-regular hypergraph of girth $g$.
\etm

Theorem \ref{t:existence_Furedi_etal} is actually stated in \cite{FLSUW95} for the  so-called \emph{bipartite biregular $(k,t,g)$-graphs}, which are bipartite graphs of girth $2g$ whose vertices in one of the partite sets have degree $k$ and degree $t$ in the other. Such a graph is precisely the incidence structure of a $k$-uniform $t$-regular hypergraph of girth $g$, see \cite{AJRS22, ErTu23} for more information on these graphs and their duality with hypergraphs. Concerning the order of $k$-uniform $t$-regular hypergraphs of given girth $g$, we cite the following result.

\btm[\cite{ElLin14}]   \label{t:d-r}
For integers $t \geq 2$ and $k, g \geq 3$, there exists a\/ $t$-regular\/ $k$-uniform linear hypergraph\/ $\cH_{t,k,g}$ of girth at least\/ $g$ on fewer than\/ $4((t-1)(k-1))^{g+1}$ vertices.
\etm

We refer to \cite{ElLin14, FLSUW95} for further information on $k$-uniform $t$-regular hypergraphs with large girth.

\subsection{Our contribution and structure of the paper}

In Section \ref{s:ind_tur} we prove that, for $t, k \geq 3$, $\turind{C_k}{K_{1,t}} \leq \frac{1}{2}(k-2)(t-1)n$, and we show that this is sharp for infinitely many values of $n$. We shall use the existence of graphs with constant link $tH$ without cycles in a certain range of lengths as a black-box, with which it is dealt further on in the paper. Doing so, we aim to make the connection between link graphs and induced Turán numbers clear and evident from the start, while focusing, from there on, to solely the constructions of various link graphs.

In Section \ref{s:existence_thms} we prove the mentioned black-box, namely, given a link graph $H$ and a graph $G$ with circumference $c(G)$ (the length of a longest cycle) and constant link $H$, then, for arbitrary integers $h \geq c(G) + 2$ and $t \geq 2$, there are infinitely many connected graphs $F$ with constant link $tH$ and no cycle of length $s$, $c(G) < s < h$. We supply two completely distinct proofs of this result. The first one uses a theorem on the existence of $t$-regular, $k$-uniform hypergraphs of given girth, combined with a coloring trick. The second proof applies what we call the ``separation lemma'' (inspired by \cite{Bla23, BlaEt23}), and a gluing technique, offering a wide range of applications.

In Section \ref{s:hyper}, we give general recursive constructions via $t$-regular $k$-uniform linear hypergraphs that produce graphs with constant link $tK_k$ and no cycles in a certain length range. We apply this technique to construct, by means of Steiner Systems, further graphs with $tK_k$-link, $k = 3, 4$, and restrictions on cycle lengths.

Section \ref{s:connected_graphs} is devoted solely to switching, an elementary operation that produces connected graphs with constant link $tL$, for any $t \geq 2$, from a given graph with constant link $L$, preserving also cycle-length restrictions, thus widening substantially the arsenal of such constructions. 

\subsection{Definitions and notation}

Our notation follows \cite{BookWest}. All the graphs in this paper are simple graphs. Given a graph $G$, we denote with $V(G)$ and $E(G)$ the set of vertices and the set of edges of $G$, respectively, while $n(G)$ is the order and $e(G)$ the number of edges of $G$. For a vertex $v \in V(G)$, we denote with $N_G(v)$ (or $N(v)$, if the  graph in question is clear from the context) the neighborhood of $v$, and with $\deg_G(v)$ ($\deg(v)$) its degree. The minimum and maximum degree of $G$ are denoted respectively with $\delta(G)$ and $\Delta(G)$. The \emph{girth} $g(G)$ and \emph{circumference} $c(G)$ are the lengths of a shortest and, respectively, a longest cycle in $G$ (if any). For a set of vertices $S \subseteq V(G)$, $G[S]$ stands for the graph induced by the vertices of $S$. Given two graphs $G$ and $G'$, $G \cup G'$ represents the disjoint union of them, while we write $tG$ for the disjoint union of $t$ copies of $G$ and $\overline{G}$ for the complement of $G$. For a vertex $v \in V(G)$, $G-v$ is the graph that is left after deleting vertex $v$ and all edges of $G$ that are incident with $v$. For integers $n, m \geq 1$, we denote with $K_n, K_{n,m}, C_n, P_n$ the complete graph on $n$ vertices, the complete bipartite graph with partite sets of sizes $n$ and $m$, the cycle of length $n$ and the path on $n$ vertices, respectively. 

Considering a hypergraph $H$, we use again $V(H)$ and $E(H)$ for its vertex and (hyper-)edge sets. For an integer $k \geq 1$, we say that $H$ is $k$-uniform if all edges are of the same size $k$. An edge of size $1$ is called \emph{trivial}. Given a vertex $v \in V(H)$, we denote with $\deg_H(v)$ ($\deg(v)$) the degree of $v$ in $H$, i.e.\ the number of hyperedges containing $v$. As above, $\delta(H)$ denotes the minimum degree of $H$. The hypergraph $H$ is $t$-regular if all its vertices have degree $t$. All cycles we consider here are those known in the literature as \emph{Berge cycles}. More precisely, a \emph{cycle} $C$ of length $n$ in $H$ is a sequence $(v_0, e_0, v_1, e_1, \ldots, v_{n-1}, e_{n-1})$ of distinct vertices $v_i \in V(H)$ and edges $e_i \in E(H)$ such that $v_i \in e_{i-1} \cap e_i$, $0 \leq i \leq n-1$, where the indices are taken modulo $n$. As above, the girth $g(H)$ of $H$ is the length of a smallest cycle of $H$ (if any). For a vertex set $S \subseteq V(H)$, we write also here $H[S]$ for the hypergraph induced by the vertices of $S$. Finally, for an edge $e \in E(H)$, we denote with $H - e$ the hypergraph that is obtained after deleting the edge $e$ from $H$, meaning that we are left with the hypergraph with vertex set $V(H)$ and edge set $E(H) \setminus \{e\}$. If $H$ is connected and $H - e$ is not connected anymore, then $e$ is called a \emph{cut-edge}.

\section{Turán numbers for induced stars and arbitrary graphs}\label{s:ind_tur}

In this section, we will show that $\turind{G}{K_{1,t}}$ is a linear function in $n$ when $G$ is not a union of $K_2$'s, while in the latter case it does not depend on $n$. That is, $\turind{G}{K_{1,t}} = \Theta(n)$, for $G \not \cong sK_2$ and any $s \ge 1$, and $\turind{sK_2}{K_{1,t}} = \Theta(1)$ for every fixed $s$ and $t$. Moreover, we will derive an explicit upper bound when $G$ is a cycle, and prove that it is sharp for an infinite number of $n$'s.
The extremal graphs will be constant-link graphs that avoid cycles of the concerned length. The existence and constructions of such graphs are precisely the matter of Sections \ref{s:existence_thms} and \ref{s:hyper}.

For the following theorem, we have to mention some concepts and related results. 
Let $G$ be a graph. We denote with $D(G)$ the family of all graphs $G - v$, for $v \in V(G)$. Observe that, if we consider $D(G)$ as a multiset, then it corresponds to the so-called deck of $G$ (see \cite{LaSc16} for a recent reference). However, for our purposes, $D(G)$ will be a simple set.  Moreover, $\chi'(G)$ denotes the chromatic index of $G$, that is, the minimum number of colors with which the edges of $G$ can be colored such that there are no two incident edges of the same color. It is well-known by Vizing's theorem \cite{Viz65} that $\Delta(G) \leq \chi'(G) \leq \Delta(G) + 1$. Finally, for graphs $G, G'$, we denote with $R(G,G')$ the \emph{Ramsey number} of $G$ and $G'$, that is the minimum integer $n$ for which every graph on at least $n$ vertices contains either a copy of $G$ or a copy of $\overline{G'}$. We will use the value of the Ramsey number $R(sK_2, K_{t})= t+2s-2$, which was determined in \cite{FaScSh80}. For a comprehensive survey about Ramsey-style theorems, we refer to \cite{Rad_survey} and the many surveys listed there as well.

\btm \label{t:indTuran-bounds}
Let $t \geq 2$ be an integer and let $G$ be any graph with $\delta(G) \ge 1$ which is not isomorphic to $sK_2$, for any $s \ge 2$. Then  
\[\left\lfloor \frac{n}{2}\right\rfloor \leq  \turind{G}{K_{1,t}} \leq \frac{1}{2}\left(R(D(G),K_{t}) – 1\right)n\]
and
\[\binom{2s-1}{2} \leq \turind{sK_2}{K_{1,t}} \leq s (t+2s-2).\]
\etm

\bpf
Assume first that $G \not\cong sK_2$, for any $s \ge 1$. Suppose $H$ is a graph of order $n$ without induced $K_{1,t}$ and containing no copy of $G$. Suppose that $\Delta(H) \ge  R(D(G), K_{t})$, and let $v$ be a vertex of maximum degree in $H$. Then, in $H[N(v)]$, there is either a member of $D(G)$ or an independent set of size $t$. Hence, there is either a copy of $G$  or an induced $K_{1,t}$ in $H[N[v]]$, a contradiction to the choice of $H$. Thus, $\Delta(H) \leq  R(D(G), K_{t}) -1$, and we obtain that
\[\turind{G}{K_{1,t}} \leq \frac{1}{2} (R(D(G)), K_{t}) -1)n.\]
For the lower bound, consider the graph $\lfloor \frac{n}{2}\rfloor K_2$ which is free from any graph $G$ with $\delta(G) \ge 1$ that is not a union of $K_2$'s.

Now suppose that $G \cong sK_2$ for some $s \geq 2$. Let $H$ be an $sK_2$-free graph of order $n$ without induced $K_{1,t}$, and let $v \in V(H)$ be a vertex of maximum degree $\Delta$. Then $\Delta \leq R(sK_2, K_{t})-1 = t+2s-3$, as otherwise the neighborhood of $v$ would contain either a copy of $sK_2$ or an independent set of size $t$, which together with $v$ would form an induced $K_{1,t}$.  Since each color class in such a coloring is a set of pairwise non-incident edges, there has to be one class of size at least $\frac{e(H)}{\chi'(H)}$, which cannot be larger than $s-1$ because $H$ is $sK_2$-free. It follows by Vizing's theorem that $e(H) \leq (s-1) \chi'(H) \leq s(\Delta +1) \leq s(t+2s-2)$. As $H$ was an arbitrary $sK_2$-free graph of order $n$ without induced $K_{1,t}$, it follows that
\[\turind{sK_2}{K_{1,t}} \leq s(t+2s-2).\]
For the lower bound, one can take $K_{2s-1}$ together with $n-2s+1$ isolated vertices. Then no induced $K_{1,2}$ occurs, and the matching number is only $s-1$.
\epf

For a graph $G$, if there is an extremal graph $H$ of order $n$ for the bound 
\begin{equation}\label{eq:ind-turan-bound}
    \turind{G}{K_{1,t}} \leq \frac{1}{2}\left(R(D(G),K_{t}) – 1\right)n
\end{equation}
of Theorem \ref{t:indTuran-bounds}, then $H$ has to be $G$-free, it does not contain induced $K_{1,t}$'s, and it has precisely $\frac{1}{2}\left(R(D(G),K_{t}) – 1\right)n$ edges. Following the proof above, we obtain that $H$ necessarily has to be regular of degree $R(D(G),K_{t}) – 1$ and such that the graph induced by the neighborhood of every vertex is isomorphic to a critical Ramsey graph on $R(D(G),K_{t}) – 1$ vertices, that is, a graph free from any member from $D(G)$ and from $\overline{K_{t}}$. Observe that this condition on a graph is necessary to be an extremal candidate but may not be sufficient as there may be copies of $G$ not fully contained in the closed neighborhood of a vertex. However, it is sufficient when $G$ has a dominating vertex, that is, a vertex of degree $n(G)-1$. Moreover, if the neighborhood of every vertex in $H$ induces a graph $L$, then $H$ is a graph with constant link $L$. It happens some times that a critical Ramsey graph is unique. Assuming this is the case for $R(D(G),K_{t})$ and $L$ is the only critical example, then bound (\ref{eq:ind-turan-bound}) can only be attained if $L$ is a link graph. If $H$ is a graph with constant link $L$ not containing $G$, then, for infinitely many values of $n \equiv 0 \pmod{n(H)}$, the upper bound is attained (just take the appropriate number of disjoint copies of $H$). That is precisely the reason why we are interested in finding certain graphs with constant link that are free of, for example, cycles of a prescribed length.

For the following theorem, in which we present a sharp upper bound on $\turind{C_k}{K_{1,t}}$, we will make use of Chvátal's theorem \cite{Chv77}, which states that $R(T,K_{t}) = (k-1)(t-1)+1$, for any tree $T$ of order $k$ and positive integers $k, t$.

\btm \label{t:tur_ind_Ck_star}
For integers $t \geq 3$ and $k \geq  3$,  
\[\turind{C_k}{K_{1,t}} \leq \frac{1}{2}(k-2)(t-1)n,\] 
and equality holds for infinitely many values of $n$.
\etm

\bpf
Observe first that $D(C_k) = P_{k-1}$ and that $R(P_{k-1}, K_{t}) = (k-2)(t-1)+1$, by Chvatal's Theorem. Hence, Theorem \ref{t:indTuran-bounds} implies that any graph $G$ on $n$ vertices that is $C_k$-free and has no induced $K_{1,t}$ must have maximum degree $\Delta(G) \leq  R(P_{k-1}, K_{t})-1 = (k-2)(t-1)$. Hence, $e(G) \leq \frac{1}{2}(k-2)(t-1)n$, proving the upper-bound.  
 
For the lower bound, let $g > k$ be an integer. Clearly, $K_{k-1}$ is a graph with constant link $K_{k-2}$ and with circumference $c = k-1$. By Theorems \ref{t:existence_const_link-gr} and \ref{t:sep_perm_constr} below, there is a graph $H$ with constant link $(t-1)K_{k-2}$ such that it is free of cycles of length $\ell$, for any $k-1 < \ell < g$. Moreover, $G$ is $(t-1)(k-2)$-regular, it has no induced $K_{1,t}$ and it has $e(G) = \frac{1}{2}(k-2)(t-1)n$ edges, where $n$ is the order of $G$. Then, for any $q \ge 1$, the graph $H_q = qG$ is $(t-1)(k-2)$-regular, it has no induced $K_{1,t}$ and it has $e(G) = \frac{1}{2}(k-2)(t-1)qn(H)$ edges. Hence, for all values of $n$ such that $n \equiv 0 \pmod{n(H)}$, there is a graph with the desired properties that attains the bound.
\epf

Let $\mathcal{C}_{k,j}$ be the family of all the graphs obtained from a cycle $C_k$ and further $j$ edges that are all incident to one particular vertex of the cycle, where $0 \leq j \leq k -2$.  
 
\begin{cor}
    Suppose $H$ is a member of $\mathcal{C}_{k,j}$. Then $\turind{H}{K_{1,t}} \leq \frac{1}{2}(k-2)(t-1)n$ and equality holds for infinitely many values of $n$.
\end{cor}

\bpf
Since $P_{k-1} \in D(H)$, Theorem \ref{t:indTuran-bounds} gives $\turind{H}{K_{1,t}} \leq \frac{1}{2}(k-2)(t-1)n$. As $H$ contains a cycle $C_k$, we necessarily have  $\turind{H}{K_{1,t}} \geq  \turind{C_k}{K_{1,t}}$. Hence, Theorem \ref{t:tur_ind_Ck_star} yields the result. 
\epf

\section{Existence of graphs with constant link $tL$ without cycles of lengths in a certain range}\label{s:existence_thms}

In this section, we will show that, provided there is a link graph $L$, the existence of connected graphs with constant link $tL$ and without cycles of lengths in a certain range is granted for any $t \geq 2$. In the first part, we show the existence without providing explicit examples, and set a relation of graphs with constant link $tL$ and hypergraphs that we will be using constantly. In the second part, we present a more explicit construction whose final gluing depends on having certain permutations, whose existence we can guarantee by a probabilistic argument.
\subsection{Existence and relation to hypergraphs}

We recall that the existence of $t$-regular, $k$-uniform linear hypergraphs with girth at least $g$, for any $t \geq 2$ and $k, g \geq 3$ is ensured by Theorem \ref{t:existence_Furedi_etal}.

Note that the girth condition already grants us with a linear hypegraph. However, as linearity will be essential in many of our proofs, we want to emphasize this property even if it is redundantly stated. Observe also that we can assume that the hypergraph given by Theorem \ref{t:existence_Furedi_etal} is connected, as otherwise we simply take one of its connected components.

By means of Theorem \ref{t:existence_Furedi_etal}, we can now give a first proof of the following theorem concerning a link graph $L$ and the existence of graphs with constant link $tL$, $t \ge 2$, which are free of cycles of a certain range of lengths.

\btm   \label{t:existence_const_link-gr}
Let\/ $F$ be a connected graph with constant link\/ $L$, order $k$ and circumference\/ $c$, where $3 \leq c \leq n$. Then, for all integers $t\geq 2$ and $g > c$, there exists a connected graph\/ $G_{t, F, g}$ with constant link $tL$ which is free of cycles of length $l$, for all $c < l < g$.
\etm

\bpf
Consider a $t$-regular $k$-uniform connected hypergraph $H$ of girth at least $g>c\geq 3$ whose existence is granted by Theorem \ref{t:existence_Furedi_etal}. Construct first a graph $G$ by taking $H$ and replacing every hyperedge with a copy of $K_k$ provided with a red-green coloring on its edges in which the green edges induce the graph $F$. This is possible due to the linearity of $H$.

Observe that $G$ has no cycles of length $\ell$ for any $k < \ell < g$. This is because if there is a cycle $C$ of length $\ell \ge k $, it necessarily contains vertices from at least two different hyperedges. So if we start at a vertex from $C$ and walk all along the cycle, we will have to cross in total at least $g$ hyperedges.

Now we delete all red edges in $G$ and name $G'$ the graph obtained in this way. As $H$ is $t$-regular and each hyperedge yields a copy of $F$ that has constant link $L$, $G'$ has constant link $tL$; here . Also, $G'$ is connected because $F$ and $H$ are connected. Since $G'$ is a subgraph of $G$, it cannot contain cycles of length $\ell$ for any $k < \ell < g$. Moreover, if it has a cycle of length $\ell \leq k$, then it has to be contained in one of the copies of $F$ whose vertices are all those of a hyperedge of $H$. It follows that $\ell \le c$. Hence, the graph $G'$ is free of $\ell$-cycles for any $c < \ell < g$. 
\epf

Provided one knows a graph $F$ of order $k$ and constant link $L$, the proof of Theorem \ref{t:existence_const_link-gr} yields also that every graph with constant link $tK_{k-1}$ (or, if wanted, every $t$-regular $k$-uniform linear hypergraph) can be used to create a graph with constant link $tL$. Moreover, information about the type of cycles contained in the different graphs involved in this construction can tell us which cycles can be avoided. Finally, one can guarantee a connected graph if the different pieces we use for this construction are connected.
\begin{prop}\label{p:coloring}
    Let $t \geq 1$, and $k \geq 2$ be integers. Let $F$ be a graph of order $k$ and constant link $L$, and let $G$ be a graph of constant link $t K_{k-1}$. Then we have the following.
    \begin{enumerate}
        \item[(i)] $G$ can be transformed into a graph $G'$ of constant link $tL$ by coloring the edges of every copy of $K_k$ in $G$ with red and green such that the green edges induce a graph isomorphic to $F$ and then deleting all the red edges.
        \item[(ii)] If $k < g$, $G$ has no cycles of length $\ell$ for any $k < \ell < g$, and $F$ has circumference $c$, then $G'$ has no cycles of length $\ell$ for any $c < \ell < g$.
        \item[(iii)] If $F$ and $G$ are connected, then $G'$ is connected as well. 
    \end{enumerate} 
\end{prop}
Observe that, for the case $t = 1$, $G \cong qK_k$ and $G' \cong qF$, for some integer $q \geq 1$, as $K_k$ is the only connected graph with constant link $K_{k-1}$. The technique used in the proof of Theorem \ref{t:existence_const_link-gr} suggests also that $t$-regular $k$-uniform linear hypergraphs and graphs with constant link $tK_{k-1}$ may be essentially the same object. This, and more, will be subject of our next theorem

It is well known that every hypergraph $H$ can be represented by its incidence graph $\iota(H)$, that is, the bipartite graph with vertex set $V(H) \cup E(H)$ and edge set $\{ve \; |\; e \in E(H), v \in e
\}$. On the other hand, every bipartite graph is the incidence relation of a hypergraph, where multiple hyperedges are allowed. Moreover, $t$-regular $k$-uniform hypergraphs of girth $g \geq 2$ correspond to biregular bipartite $(k,t,g)$-graphs, i.e. bipartite graphs of girth $2g$ whose vertices in one of the partite sets have degree $k$ and degree $t$ on the other.

Let $t, k, g$ be integers with $t, k \geq 2$, and $g \geq \max\{4, k+1\}$. We denote with $\mathcal{L}_{k,t}$ the set of all graphs of constant link $tK_{k-1}$, and let $\mathcal{L}_{k,t,g}$ be the subset of $\mathcal{L}_{k,t}$ containing all those graphs having cycles of length $g$ but no cycles of length $\ell$ for any $\ell < g$ except for those contained in the cliques of size $k$. Let $\mathcal{B}_{k,t}$ be the set of all biregular bipartite graphs with degrees $k$ and $t$ and girth at least $8$, and $\mathcal{B}_{k,t,g}$ the set of all those having girth $2g$. Let $\mathcal{H}_{k,t}$ be the set of all $t$-regular $k$-uniform hypergraphs of girth at least $4$, and more restrictively $\mathcal{H}_{k,t,g}$ the set of those having girth $g$. 

We show in the next result that the sets $\mathcal{B}_{k,t}$, $\mathcal{H}_{k,t}$ and $\mathcal{L}_{k,t}$ are essentially the same, as well as the sets $\mathcal{B}_{k,t,g}$, $\mathcal{H}_{k,t,g}$, and $\mathcal{L}_{k,t,g}$, for $t, k \geq 2$, and any $g \geq \max\{4, k+1\}$. Observe that, as the hypergraphs in $\mathcal{H}_{k,t}$ are linear, we can easily replace the hyperedges in any $H \in \mathcal{H}_{k,t}$ by copies of $K_k$ to create a graph $\varphi(H)$. Since $H$ is $t$-regular, every vertex $v$ is contained in $t$ hyperedges whose intersection is precisely $\{v\}$. Moreover, as $H$ has girth at least $4$, there are no hyperedges containing vertices across the  hyperedges incident with $v$. Hence, it is clear that, in $\varphi(H)$, the neighborhood of every vertex induces a $t K_{k-1}$. Consequently, the mappings
\vskip-1ex
  $$\iota: \mathcal{H}_{k,t} \to \mathcal{B}_{k,t}, \;\; H \mapsto \iota(H)
  $$
\vskip-2ex
\noindent  \mbox{and}
\vskip-4ex
  $$
    \varphi: \mathcal{H}_{k,t} \to \mathcal{L}_{k,t},\;\; H \mapsto \varphi(H)
    $$
are well defined, and $\iota$ is bijective. We will show that $\varphi$ is bijective as well and, for $g > k \geq 3$, the restriction of these mappings to $\mathcal{H}_{k,t,g}$ yields a one-to-one correspondence between the families $\mathcal{B}_{k,t,g}$, $\mathcal{H}_{k,t,g}$, and $\mathcal{L}_{k,t,g}$.

\btm \label{t:families_equiv}
Let $t, k, g$ be positive integers, where $t, k \geq 2$, and $g \geq \max\{4, k+1\}$. Then there is a one-to-one correspondence between the families $\mathcal{B}_{k,t}$, $\mathcal{H}_{k,t}$ and $\mathcal{L}_{k,t}$,
and also between the families $\mathcal{B}_{k,t,g}$, $\mathcal{H}_{k,t,g}$ and $\mathcal{L}_{k,t,g}$.
\etm

\bpf
That there is a one-to-one correspondence between the families $\mathcal{H}_{k,t}$ and $\mathcal{B}_{k,t}$, and between the families $\mathcal{H}_{k,t,g}$ and $\mathcal{B}_{k,t,g}$ is clear by means of the incidence structure explained above, corresponding to the mapping $\iota$ and its restriction to $\mathcal{H}_{k,t,g}$. We will prove now that also $\varphi$ is bijective. Let $G \in \mathcal{L}_{k,l}$. We show next that every clique of size $k$ can be replaced by a hyperedge such that a hypergraph in $\mathcal{H}_{k,t}$ is created, thus obtaining a reverse operation to $\varphi$. To this aim, we need to prove that every pair of cliques $S_1$, $S_2$ of size $k$ can intersect in at most one vertex. Suppose to the contrary that $x, y \in S_1 \cap S_2$, where $x \neq y$. Then the set $(S_1 \cup S_2) \setminus \{x\}$, which is a subset of $N_G(y)$, induces a connected graph on at least $k$ vertices (as the cliques are different). But the neighborhood of any vertex is isomorphic to $tK_{k-1}$, whose connected components have $k-1$ vertices, a contradiction. Hence, $|S_1 \cap S_2| \leq 1$. Now it is easy to see that every copy of $K_k$ can be replaced by a hyperedge. Since the closed neighborhood of every vertex $v$ is a collection of $t$ $K_k$'s that intersect in $v$, it is clear that there are no more edges left and that this way a $t$-regular $k$-uniform hypergraph is created. Finally, observe that the girth of such a hypergraph is at least $4$. This is due to the fact that, for any two cliques $S_1$, $S_2$ from $G$ that intersect in one vertex $x$, there are no edges between the sets $S_1 \setminus \{x\}$ and $S_2 \setminus \{x\}$. Hence, $\varphi: \mathcal{H}_{k,t} \to \mathcal{L}_{k,t}$ is a bijective mapping.

Let $g > k$. We know already that $\iota(\mathcal{H}_{k,t,g}) = \mathcal{B}_{k,t,g}$. Now we will show that $\varphi(\mathcal{H}_{k,t,g}) = \mathcal{L}_{k,t,g}$. Let $H \in \mathcal{H}_{k,t, g}$. Observe first that any cycle of length $\ell$ in $H$ translates to a cycle of length $\ell$ in $\varphi(H)$. We will show that $\varphi(H)$ has no cycles of length $\ell < g$ other than those contained in the cliques of size $k$. Suppose there is a cycle $C$ in $\varphi(H)$ of length $\ell < g$ not all of its vertices belonging to a clique of size $k$. It may happen that more than two consecutive vertices in the cycle come from the same hyperedge in $H$. But, as they are not all contained in a clique of size $k$, there must be several hyperedges involved in this cycle (not just one). These hyperedges form a cycle in $H$, which has length at least $g$. Thus, the cycle $C$ has length at least $g$ and we have shown that $\varphi(\mathcal{H}_{k,t,g}) \subseteq \mathcal{L}_{k,t,g}$. Now consider a graph $G \in \mathcal{L}_{k,t,g}$ and its preimage $H = \varphi^{-1}(G) \in \mathcal{H}_{k,t}$. Again, cycles that are not contained in the cliques of size $k$ in $G$ are mapped to cycles in $H$. Since, by means of $\varphi^{-1}$, we are replacing every clique of size $k$ with a hyperedge, and as every cycle of $G$ of length at most $k$ is contained a clique of size $k$, all cycles in $H$ have length at least $g$. Clearly, the presence of a cycle of length at least $g$ in $G$ yields a cycle of length at least $g$ in $H$ (although the former may be much longer than the latter). Thus $H \in \mathcal{H}_{k,t,g}$, and we have shown that $\varphi^{-1}(\mathcal{L}_{k,t,g}) \subseteq \mathcal{H}_{k,t,g}$. Altogether, we have that $\varphi(\mathcal{H}_{k,t,g}) = \mathcal{L}_{k,t,g}$, yielding finally the one-to-one correspondence between the families $\mathcal{B}_{k,t,g}$, $\mathcal{H}_{k,t,g}$ and $\mathcal{L}_{k,t,g}$.
\epf

\subsection{Construction via separating permutations}\label{s:permut}

In the following, we will present an interesting construction that, given integers $k$ and $g$ such that $k < g$, yields a graph with constant link $tK_k$ avoiding cycles of lengths $h$ for all $k < h < g$. Our method is to take a combination of explicitly constructed parts and a probabilistic lemma that deals with the existence of particular \emph{separating permutations}. 

To be able to present the lemma, we first need some terminology. Let $X$ be a set of $n$ elements. By a cyclic permutation over $X$ we mean a bijective mapping $\pi$ from $X$ to the vertex set of the $n$-cycle $C_n=v_1v_2\dots v_n$. If $\pi(x)=v_i$ and $\pi(y)=v_j$, then the distance $\dist_\pi(x,y)$ of $x$ and $y$ under $\pi$ is meant as the graph-theoretic distance, $\min(|i-j|,n-|i-j|)$ between the vertices $v_i$ and $v_j$ along $C_n$.

\begin{lemma} \label{l:sep_perm}
Let $d\geq 2$ be an integer and $X$ a set of $n$ elements. If $n > 10 (d+1)^8$, then there exist $d$ cyclic permutations $\pi_1,\dots,\pi_d$ over $X$ such that, for any two $x,y\in X$, $\dist_{\pi_\ell}(x,y)\leq d$ is valid for at most one $\pi_\ell$ ($1\leq\ell\leq d$).
\end{lemma}

\bpf
Let $N:=n-1+10d^4$, and let $X'$ be
 a set of $N+1=n+10d^4$ elements containing $X$.
We make a random selection of $d$ independent cyclic
 permutations $\pi'_1,\dots,\pi'_d$ over $X'$, in the sense of
 total independence; i.e., the choice of each
 $\pi'_i$ has no correlation with any combination of the
 choices of the other $d-1$ permutations, and
 each of the cyclic permutations $\pi'$ over $X'$ is equally
 likely to coincide with $\pi'_i$.

At first, we put a requirement stronger than the one in the lemma.
Call an element $x\in X'$ bad if, for some $y\in X'\smin\{x\}$,
 there are at least two indices $i,j$ such that
  $\dist_{\pi'_i}(x,y)\leq d+2$ and
 $\dist_{\pi'_j}(x,y)\leq d+2$.
We claim that the probability of an $x$ to be bad is small.
Indeed, the two indices $i,j$ making $x$ bad can be chosen in
 at most $\binom{d}{2}$ ways, and there are $N=n-1+10d^4$
 candidates for the element $y$ too close to $x$.
This $y$ should be placed as one of the $2d+4$ elements at
 distance at most $d+2$ apart from $x$ in both $\pi'_i$ and $\pi'_j$.
Hence, the probability that $x$ is bad is generously bounded above by
\begin{equation}\label{eq:p}
    p := \ppp(x \mathrm{ ~is~bad }) \leq
   N\cdot \binom{d}{2} \cdot \left(\frac{2d+4}{N}\right)^2 
   \leq \frac{4d^4}{N}
    < \frac{5d^4}{N+1} \,,
\end{equation}
  moreover the expected number of bad elements is $(N+1)p$, thus
\[ \eee(|x\in X' : x \mathrm{ ~is~bad }|) < 5d^4.\]
Consequently, by Markov's inequality we obtain
  \begin{equation}   \label{eq:few-x}
  \ppp(|x\in X' : x \mathrm{ ~is~bad }| < 10 d^4) > 1/2 \,.
  \end{equation}

Consider next a permutation $\pi'_i$, for some $i \in [d]$. For $j \in [N+1]$, let $S_{i,j}$ be the set of $d+4$ consecutive elements in the cyclic order determined by $\pi'_i$ starting at $j$.
We have $d$ possible choices for $i$ and $N+1$ choices for $j$
 (i.e., where to start the segment $S_{i,j}$).
Call an $S=S_{i,j}$ bad if it contains at least three bad elements.
There are $\binom{d+4}{3}$ ways to select three elements
 $x,y,z$ from $S$; if the number of bad vertices is $b$, then
 the probability that three bad elements are mapped to $x,y,z$
 is $\frac{b(b-1)(b-2)}{(N+1)N(N-1)}$ which is smaller than
 $p^3$; and, as told above, there are $(N+1)d$
 possible positions of $S$.
Consequently, using (\ref{eq:p}) and the fact that $N+1 > 10 (d+1)^8$, the probability that no bad segments of size
 $d+4$ occur is at least
 \begin{equation}   \label{eq:no-S}
  1 - (N+1)\,d \binom{d+4}{3} p^3 > 1 - \binom{d+4}{3} \frac{125d^{13}}{(N+1)^2} >  1 - \frac{125 (d+1)^{16}}{3 (N+1)^2}  >  7/12.
 \end{equation}

Inequalities (\ref{eq:few-x}) and (\ref{eq:no-S}) together
 imply the existence of $d$ permutations $\pi'_1,\dots,\pi'_d$
 over $X'$ such that the size $b$ of the set $B$ of bad
 elements is smaller than $10d^4$.
Since $n$ is much bigger than $10d^4\cdot(2d+5)$, we can extend
 $B$ to a set $B'$ of size exactly $10d^4=(N+1)-n$ such that each
 $S=S_{i,j}$ meets $B'$ in at most two elements.
By this condition, for any $t>0$, any $d+2+t$
 consecutive elements of $X'$ under any $\pi'_\ell$ contain
 at most $t$ elements from $B'$.

Fix a bijective mapping $\psi:X'\smin B'\to X$ and, for
 $\ell=1,\dots,d$, let $\pi_\ell=\psi(\pi'_\ell|_{X'\smin B'})$,
 where
 $\pi'_\ell|_{X'\smin B'}$ is obtained from $\pi'_\ell$
 by removing the elements of $B'$ from the cyclic order of
 $X'$ determined by $\pi'_\ell$ (while keeping the order of
 the remaining elements unchanged).
As shown above, if two elements of $X'\smin B'$ have
 distance greater than $d+2$ under $\pi'_\ell$, then they are
 mapped by $\psi$ to two elements of $X$ that have
 distance greater than $d$ under $\pi_\ell$.
Since $X'\smin B'$ does not contain any bad elements,
 we obtain that the permutations $\pi_1,\dots,\pi_d$
 satisfy the requirements of the lemma.
\epf

Given integers $k\geq 2$ and $t\geq 2$, we provide now a construction of graphs with constant link $tK_{k-1}$. In the final step of this construction, we need to glue together a series of graphs. We will show in Theorem \ref{t:sep_perm_constr} that this gluing can be done avoiding cycles in a certain range of lengths. To this aim, we define the following product of graphs. Let $G$ and $H$ be graphs and let $S \subseteq V(G)$ be a set of vertices in $G$. We define $G \vee_S H$ as the graph obtained from a copy of $G$ and $|S|$ copies of $H$, say $H_v$, $v \in S$, where each vertex $v \in S$ is joined by an edge to each vertex of $H_v$.

\begin{constr} \label{const:k-ary-tree-like}
Let $t, k, q \geq 2$ be integers. Let $m_i = (k-1)^{i} (t-1)^{i-1}$, for $i \ge 1$. Let $\Pi = \{\pi_1, \pi_2, \ldots, \pi_t\}$ be a set of $t$ permutations selected from the symmetric group $S_{m_q}$ of order $m_q$. We construct a graph $G(t,k,q,\Pi)$ with constant link $tK_{k-1}$ in the following way.
\begin{enumerate}
    \item[(i)] Let $G_1 \cong K_k$ and let $x \in V(G_1)$ be a vertex called the root of $G_1$. Let $W_1 = V(G) \setminus \{x\} = \{\binom{1}{\ell} \;:\; \ell \in [k-1]\}$. That is, we label all neighbors of $x$ by means of $2 \times 1$ matrices. 
    
    \item[(ii)] Define recursively, for $i \ge 2$, $G_i = G_{i-1} \vee_{W_{i-1}} (t-1)K_{k-1}$, $V(G_i) = V(G_{i-1}) \cup W_i$, where
    \[W_i = \left\{\binom{u \; u_i}{v \; v_i}\; :\; u_i \in [t-1], v_i \in [k-1], \binom{u}{v} \in W_{i-1} \right\}.\]
    The assignation of names of the new vertices by means of $2 \times i$ matrices is done in such a way that the matrix $\binom{u \; u_i}{v \; v_i}$ represents the $v_i$-th neighbor in the $u_i$-th copy of $K_{k-1}$ of the vertex $\binom{u}{v} \in W_{i-1}$. In this way, the vertices in $W_i$ are linearly ordered by the rule $\binom{u_1 u_2 \ldots u_i}{v_1 v_2 \ldots v_i} \leq \binom{u'_1 u'_2 \ldots u'_i}{v'_1 v'_2 \ldots v'_i}$ if and only if $u_{\ell} \leq u'_{\ell}$ and $v_{\ell} \leq v'_{\ell}$ for all $1 \leq \ell \leq i$.  For $1 \leq \ell \leq i$, we will say that the elements of $W_\ell$ are the vertices in level $\ell$ of $G_i$. Observe that $|W_i| = m_i$, for any $i \ge 1$. 
    
    \item[(iii)] Take $t$ copies of $G_q$, say $G_q^1, G_q^2, \ldots, G_q^t$ with roots $x_1, x_2, \ldots, x_t$, and let $W_q^j$ be the set of vertices in level $q$ of $G_q^j$, $1 \le j \le t$. The vertices in $W_q^j$ have the labels $\binom{u}{v}$, with $u \in [t-1]^q$ and $v \in [k-1]^q$ according to the linear order given above. Let $\varphi_j : W_q^j \to [m_q]$ be a bijective mapping following that order, that is, $\binom{u}{v} \leq \binom{u'}{v'}$ if and only if $\varphi_j(\binom{u}{v}) \leq \varphi_i(\binom{u'}{v'})$.
    
    \item[(iv)] Identify all roots  $x_1, x_2, \ldots, x_t$ to one vertex $x$.
    
    \item[(v)] For every $j \in [m_q]$, identify the vertices $\varphi_i^{-1}(\pi_i(j))$, $1 \leq i \leq t$.
    \end{enumerate}
\end{constr}

\begin{figure}[h]
    \centering
    \includegraphics[width=0.9\linewidth]{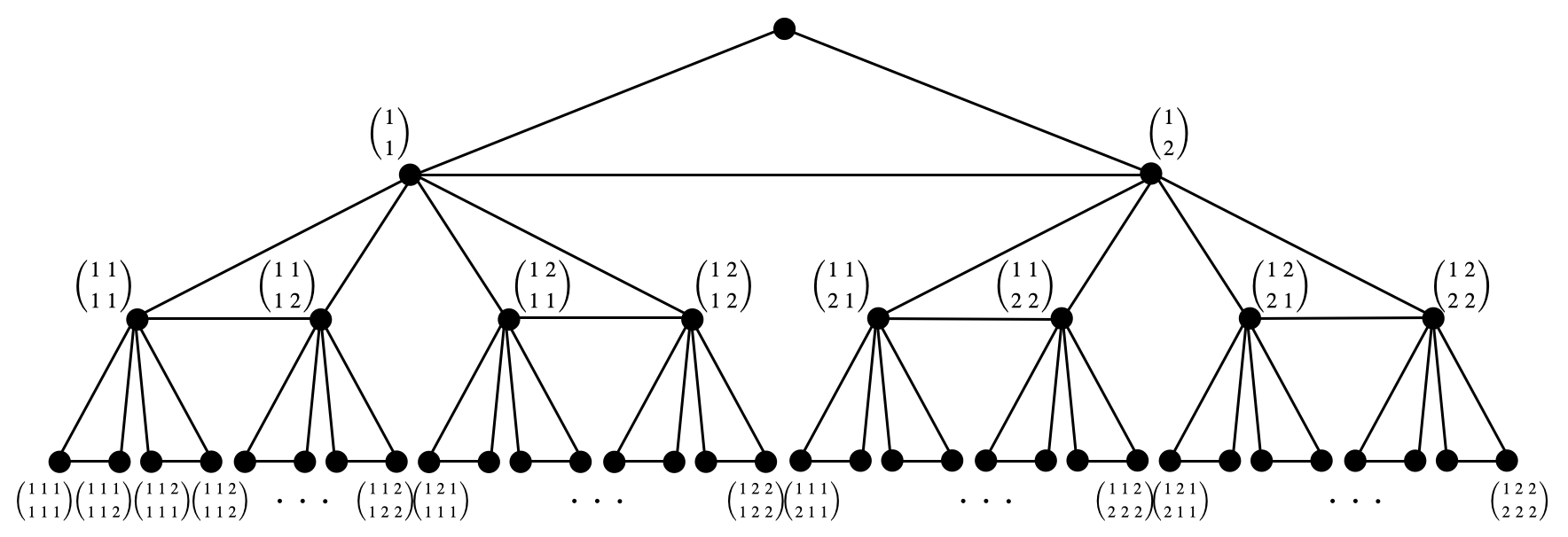}
    \caption{Example of labeling of $G_3$ in Construction \ref{const:k-ary-tree-like}, case $t = 3$, $k = 3$.}
    \label{f:labeling_constr10}
\end{figure}

\begin{figure}[h]
    \centering
    \includegraphics[width=0.6\linewidth]{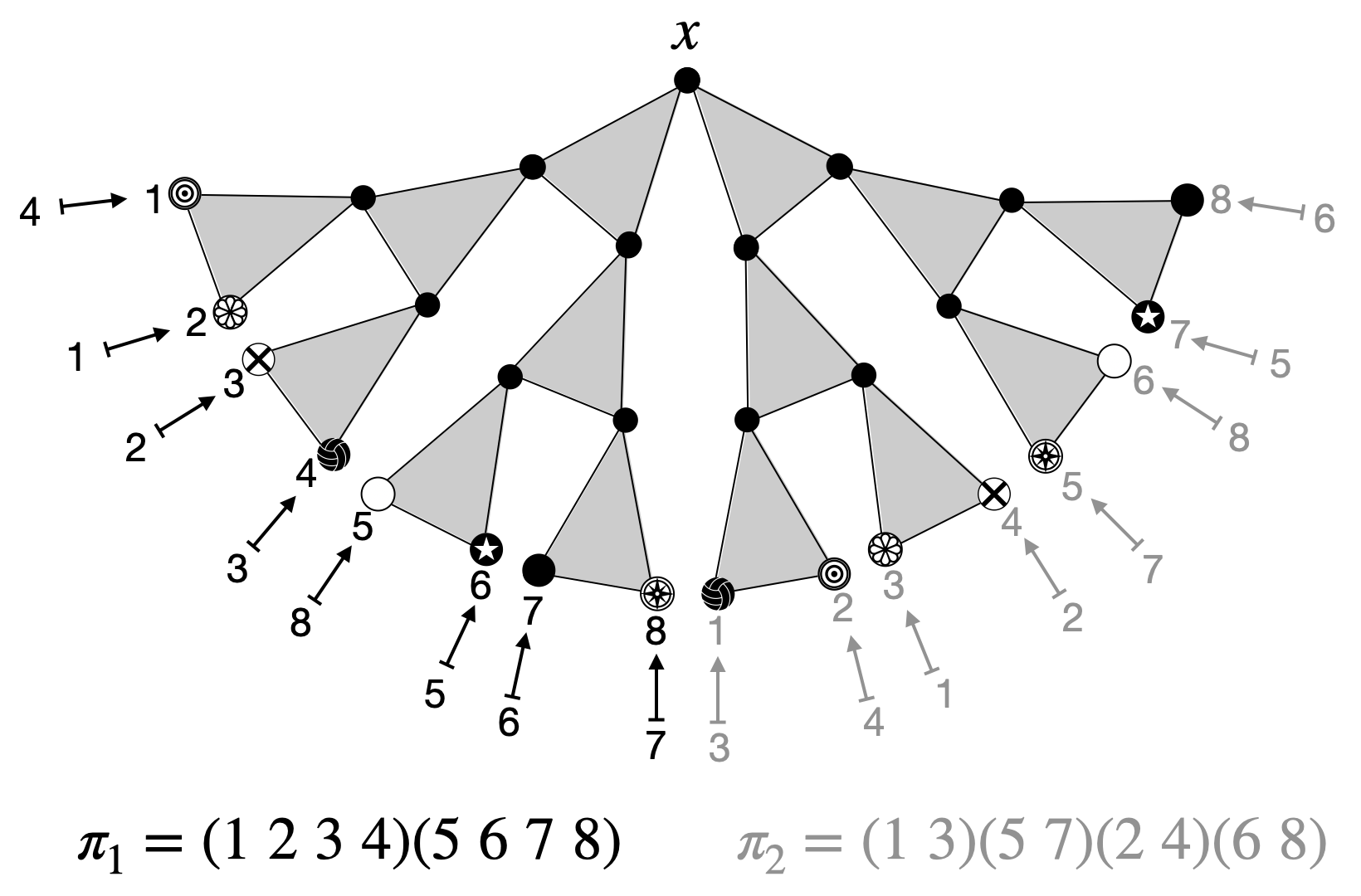}
    \caption{Example of gluing of two copies of $G_3$ in Construction \ref{const:k-ary-tree-like}, case $t = 2$, $k = 3$. Vertices in the last layer with the same shape are identified.}
    \label{f:gluing_constr10}
\end{figure}

See Figures \ref{f:labeling_constr10} and \ref{f:gluing_constr10} for examples concerning the labeling and gluing in Construction \ref{const:k-ary-tree-like}. We prove in the next theorem that there can be chosen a large enough integer $q$ and a set $\Pi \subseteq S_{m_q}$ of permutations such that the graph $G(t,k,q,\Pi)$ has no cycles in a certain range of lengths.

\btm \label{t:sep_perm_constr}
Let $t, k, g$ be integers such that $t, g > k \geq 2$. Then there is a large enough integer $q$ and a set $\Pi = \{\pi_1, \pi_2, \ldots, \pi_t\}$ of $t$ permutations in $S_{m_q}$, where $m_q = (k-1)^{q} (t-1)^{q-1}$, such that the graph $G(t,k,q,\Pi)$ has no cycles of length $h$ for any $k < h < g$.
\etm

\bpf
Let $m_i = (k-1)^{i} (t-1)^{i-1}$, for $i \ge 1$. Observe that the blocks (that is, the maximal $2$-connected components) of the graph built up to step (v) of Construction \ref{const:k-ary-tree-like} are all isomorphic to $K_k$. Hence, the only cycles have been created up to this point are of length at most $k$. Any other cycle in this construction appears in step (v), that is, at the moment of gluing vertices in the last level $q$ of the different $G_q^i$'s, $1 \leq i \leq t$. In order to avoid the cycle lengths $h$, for $k < h < g$, we will force that any two vertices of level $q$ in any of the copies $G_q^i$ that have a common ancestor in level $q-g$ or higher (meaning closer), will be glued to vertices of the other $G_q^j$'s, $1 \leq j \leq t$, $j \neq i$, whose closest common ancestor is sufficiently far away, which, to be on the safe side, we will say it shall be at a level less than $q-g$. Observe that pairs of vertices in the last level of any of the $G_q^i$'s have a common ancestor in level $q-g$ or higher precisely if they are contained inside any interval $\varphi^{-1}([s,s+m_g-1])$, for an $s \in [m_q - m_g+1]$. Hence, we want to choose $t$ permutations in $S_{m_q}$ that prevent that two or more elements inside the intervals $[s,s+m_g-1]$, $s \in [m_q - m_g+1]$, are associated with elements that are too close to each other. This is indeed guaranteed by choosing $t$ permutations $\pi_1, \pi_2, \ldots, \pi_t$ such that, for any two $x,y\in [m_q]$, $\dist_{\pi_\ell}(x,y)\leq m_g$ is valid for at most one $\pi_\ell$ ($1\leq\ell\leq t$). Thus, taking $d = \max \{m_g, t\}$, and letting $q$ be large enough such that $m_q > 10 (d+1)^8$, we apply Lemma~\ref{l:sep_perm}, by which it follows that there are $d$ permutations $\pi_1, \pi_2, \ldots, \pi_d \in S_{m_q}$ such that, for any $x, y \in [m_q]$, $\dist_{\pi_i}(x,y)\leq d$ is valid for at most one $\pi_i$, $1 \leq i\leq d$. As $d \geq m_q$, it also holds that $\dist_{\pi_i}(x,y) \leq m_g$ is valid for at most one $\pi_i$, $1 \leq i\leq d$. Selecting $t$ from these permutations, we obtain the desired set $\Pi$ such that the graph $G(t,k,q,\Pi)$ has no cycles of length $h$ for any $k < h < g$.
\epf

\section{Recursive constructions of graphs with constant link $tL$ via known hypergraphs}\label{s:hyper}

We will present in this section a general recursive method to generate graphs with constant link $tL$, 
for any $t\geq 2$ and for any admissible link $L$, when a graph with constant link $L$ is available.
The method will require the existence of certain regular and uniform hypergraphs. Moreover, it will yield a graph with constant link $tL$ free from cycles of a certain range of lengths, as long as the employed hypergraphs have large enough girths. The following lemma gives the essence of the recursive construction.

\begin{constr}\label{gral_constr}
    Let $F$ be a graph with constant link $L$, and let $k = n(F) \geq 3$. Suppose $G$ is a graph of order $n$ and constant link $tL$, for some $t \ge 1$. Let $H$ be an $n$-regular $k$-uniform hypergraph of order $N$.
    We construct a graph $\Gamma = \Gamma(G,H)$ of order $nN$ and constant link $(t+1)L$ in the following way.
    \begin{enumerate}
        \item[(i)] At each vertex $x \in V(H) = \{x_1, x_2, \ldots, x_N\}$, we enumerate the hyperedges incident with $x$ as $E_{x,1},\dots, E_{x,n}$.
        \item[(ii)] Take $N$ disjoint copies of $G$ on the vertex sets $V_1, V_2, \ldots, V_N$, where $V_i = \{v_{i,1}, v_{i,2}, \ldots, v_{i,n}\}$, for $1 \leq i \leq N$, and set $V(\Gamma) = V_1 \cup V_2 \cup \cdots \cup V_N$.
        \item[(iii)] With each edge $E\in \mathcal{E}(H)$, a $k$-element subset $f(E)$ of $V(\Gamma)$ is associated by the rule $v_{i,j}\in f(E)$ if and only if $E_{x_i,j} = E$.
        \item[(iv)] For each hyperedge $E \in \mathcal{E}(H)$, place a copy of $F$ on the set $f(E)$.
    \end{enumerate}   
\end{constr}

The following theorem gives a general recursive method using Construction \ref{gral_constr} to build graphs with constant link $tL$, provided that a graph with constant link $L$ and a certain sequence of uniform regular
hypergraphs is available. If these hypergraphs have additionally sufficiently large girth, then cycles of a range of lengths can be avoided. Let us note that girth assumptions also exclude 2-cycles; i.e., hypergraphs of that kind are necessarily linear.

\btm \label{t:general_method}
Let $F$ be a graph of constant link $L$ and order $k \geq 3$. Let $(H_t)_{t \ge 1}$ be a sequence such that each $H_t$ is a $k$-uniform, $n_t$-regular hypergraph of order $v_t$, and suppose that $n_1=k$ and $n_t = k \Pi_{i = 1}^{t-1} v_i$ holds for every $t \ge 2$. Let $G_1 = F$, and define $G_{i+1} = \Gamma(G_i,H_i)$, for $i \ge 1$. Then the following holds for every $t \ge 1$.
\begin{enumerate}
    \item[(i)] The graph $G_t$ has constant link $tL$ and order $n_t$.
    \item[(ii)] If $F$ has circumference $c$, and $H_i$ has girth $g > c/2$ for $1\leq i\leq t-1$, then $G_t$ has no cycles of length $\ell$, for any $c < \ell < 2g$.
    \item[(iii)] If $G_t$ and $H_t$ are connected, then $G_{t+1}$ is connected, too. 
\end{enumerate}
\etm

\bpf
Item (iii) is easily verified. To prove (i), we will proceed by induction on $t$. For $t = 1$, it holds by assumption. Assume now $G_t$ has constant link $tL$ and order $n_t$ for some $t \ge 1$. By Construction~\ref{gral_constr}, $G_{t+1} = \Gamma(G_t,H_t)$ has constant link $(t+1)L$ because the substitutions of $F$ into the copies of edges of $H_t$ are done in a vertex-disjoint way, hence with just one additional copy of $F$ at each vertex. Moreover, $G_{t+1}$ is of order 
\[n_t n(H_t) = \left(k \prod_{i = 1}^{t-1} v_i \right) v_{t} = k \prod_{i = 1}^{t} v_i = n_{t+1}.\]
Suppose now that $H_t$ has girth $g > c/2$, for every $t \ge 1$. 
Again by induction on $t$, we will prove that $G_t$ has no cycles of length $\ell$, for any $c < \ell < 2g$. This is clearly valid for $G_1$.
Assuming that the same holds for $G_t$, we will prove that also $G_{t+1}$ has no cycles in this range of lengths. Consider $G_{t+1} = \Gamma(G_t,H_t)$ as is given in Construction \ref{gral_constr}, with its vertex set $V(G_{t+1}) = V_1 \cup V_2 \cup \cdots \cup V_{v_t}$, a copy of $G_t$ placed on every set $V_i$ and a copy of $F$ on every set $f(E)$, where $E \in \mathcal{E}(H_t)$. Recall also that the sets $f(E)$ are pairwise disjoint. Suppose there is a cycle $C\subset G_{t+1}$ of length $\ell > c$. By construction, two consecutive vertices from $C$ are contained either in the same set $V_i$ or in the same set $f(E)$ for some unique $E \in \mathcal{E}_t$ (as $f$ transforms the hyperedges of $H_t$ to disjoint ones). Since a longest cycle in $F$ is of length $c$, we know that $C$ cannot be fully contained in any set $V_i$. So we can assume that $C$ is the cyclic concatenation of some number $2q\geq 2$ of paths $P_1,P_2,\dots,P_{2q}$ such that, for $s=1,\dots,q$, $P_{2s-1}$ is inside some $V_{i_s}$ and each $P_{2s}$ is inside $f(E)$, for some $E \in \mathcal{E}(H_t)$. Let $u_{2s},u_{2s+1}$ be the first and last vertex in $P_{2s}$, and let $E_{2s} \in \mathcal{H}_t$ be the hyperedge such that $u_{2s},u_{2s+1} \in f(E_{2s})$, $1 \leq s \leq q$. Then $(E_2, E_4, \ldots, E_{2q})$ is a closed trail in $\mathcal{H}_t$. But $H_t$ has girth at least $g$, by assumption. Thus, $\ell\geq 2q \geq 2g$, as required.
\epf

\subsection{Constructions from Steiner systems}\label{s:stein}

In this section, we will use Construction \ref{gral_constr} and Theorem \ref{t:general_method} to construct graphs of constant link $tK_k$, $t \geq 1$, via Steiner systems. A \emph{Steiner system} $S(s, k, v)$ is a pair $(V, B)$ where $V$ is a set of $v$ elements called \emph{points} and $B$ is a family of $k$-element subsets of $V$ called \emph{blocks}, such that each $s$-element subset of $V$ is contained in exactly one block. It is well known \cite{Wil75} that, assuming $v$ to be sufficiently large with respect to $k$, a Steiner system $S(2,k,v)$ exists if and only if
 \begin{equation}   \label{eq:divisibility_Stei}
     k(k-1) \mid v(v-1)
     \qquad \mathrm{and}\qquad
     k-1 \mid v-1 \,.
 \end{equation}
 
\nin

Concerning $3 \leq k \leq 9$, detailed research has been carried out investigating the possible exceptions $v$ for which the necessary conditions of (\ref{eq:divisibility_Stei}) are not sufficient for the existence of Steiner systems $S(2,k,v)$;
see Table~I in \cite[p.~375]{AdBryBu07} or Table 3.3 in \cite[pp.~72--73]{CD-hand}. In particular, no exceptional values occur for $3 \leq k \leq 5$. Observe that Steiner systems $S(2,k,v)$ can be regarded as $k$-uniform, $t$-regular, linear (hence, girth-3) hypergraphs with $t = \frac{v-1}{k-1}$.  For more information on Steiner systems, which are a special family of combinatorial designs, see \cite{CD-hand, Sti04}.

We can now prove the following existence result. 

\btm\label{t:Steiner}
Let $k \geq 3$, $q\geq 1$, and $v_1,v_2,\dots,v_{t-1}$ be integers such that there exist Steiner systems $H_i=S(2, k, v_i)$
 for all $1\leq i\leq t-1$.
Denoting $n_i=\frac{v_i-1}{k-1}$, assume that $n_1=qk$ and
 $n_j = qk \Pi_{i = 1}^{j-1} v_i$ holds for all
  $2\leq j \leq t-1$.
We recursively define
\[S_1(k) = qK_k, \mbox{ and \ } S_{i+1}(k) = \Gamma(S_i(k), H_i(k))\]
 for $i=1,\dots,t-1$.
Then $S_t(k)$ is a graph of order $n_t =  k q \prod_{i=1}^{t-1}v_i$ and constant link $tK_{k-1}$. 
\etm

\bpf
As noted above, Steiner systems $S(2,k,v_i)$ are $k$-uniform
 $\frac{v_i-1}{k-1}$-regular (i.e., $n_i$-regular) hypergraphs.
Clearly, $S_1(k) = qK_k$ has constant link $K_{k-1}$. Hence, by Theorem \ref{t:general_method}, $S_t(k)$ has order $n_t = qk \prod_{i=1}^{t-1}v_i$ and constant link $tK_{k-1}$, for all $t \ge 1$. 
\epf

\begin{rmk}
Wilson's theorem implies that the Steiner systems $S(2,k,v_i)$ required
 in the conditions of Theorem \ref{t:Steiner} surely exist
 whenever $q$---or, equivalently, $v_1$---is sufficiently large.
Note that the necessary condition $(\ref{eq:divisibility_Stei})$ holds
 for all $q\geq 1$ with $v:=v_1=qk(k-1)+1$.
\end{rmk}

Observe however that, for $q > 1$, we cannot guarantee connection. 
On the other hand, for the cases that $k = 3,4,5$, we can start the recursion with just one copy of $K_k$, i.e., $q=1$ and setting $S_1(k) = K_k$, and so the connectivity of $S_t(k)$ is guaranteed. 
Further, for $k = 3, 4$, we can exclude some cycle lengths, too.

\btm\label{t:small_Steiner}
Let $k \in \{ 3, 4, 5\}$. 
Set $v_1=k^2-k+1$, $v_j = k (k-1)  \prod_{i=1}^{j-1}v_i + 1$ for $1<j<t$, and let $H_i(k)$ be an $S(2,k, v_i)$-system, for all $i \geq 1$. We define
\[S_1(k) = K_k, \mbox{ and } S_{i+1}(k) = \Gamma(S_i(k), H_i(k))\]
 for $i=1,\dots,t-1$.
Then $S_t(k)$ is a connected graph of order $n_t =  k \prod_{i=1}^{t-1}v_i$ and constant link $tK_{k-1}$, for any $t \ge 1$. Moreover, $S_t(3)$ is $\{C_4,C_5\}$-free, and $S_t(4)$ is $C_5$-free, for every $t \ge 1$.
\etm

\bpf
We know that, for $k = 3,4,5$, there is a Steiner system $S(2,k,v)$ for any $v$ such that $v \equiv 1 \pmod{k(k-1)}$. Hence, we can set $v_1 = k(k-1)+1$ in Theorem \ref{t:Steiner}, and it follows that $q = 1$. Moreover, $S_t(k)$ has order $n_t =  k \prod_{i=1}^{t-1}v_i$ and constant link $tK_{k-1}$, for all $t \geq 1$. 

As $S_1(k) = K_k$ is connected, and $H_t(k)$ is connected as well for all $t \geq 1$, it follows by item (iii) of Theorem \ref{t:general_method}, using again induction, that $S_t(k)$ is connected for all $t\ge 1$.

Since in $H_t(k)$ every pair of vertices is contained in exactly one hyperedge, we have $g(H_t(k)) = 3$. Moreover, $S_1(k)$ has circumference $k$ and $k < 2g(H_t(k)) = 6$ when $k \leq 5$. Hence, it follows by item (ii) of Theorem \ref{t:general_method} that $S_t(3)$ is $\{C_4,C_5\}$-free, and $S_t(4)$ is $C_5$-free, for every $t \ge 1$.
\epf

\begin{figure}[ht]
    \centering
    \includegraphics[width=0.7\linewidth]{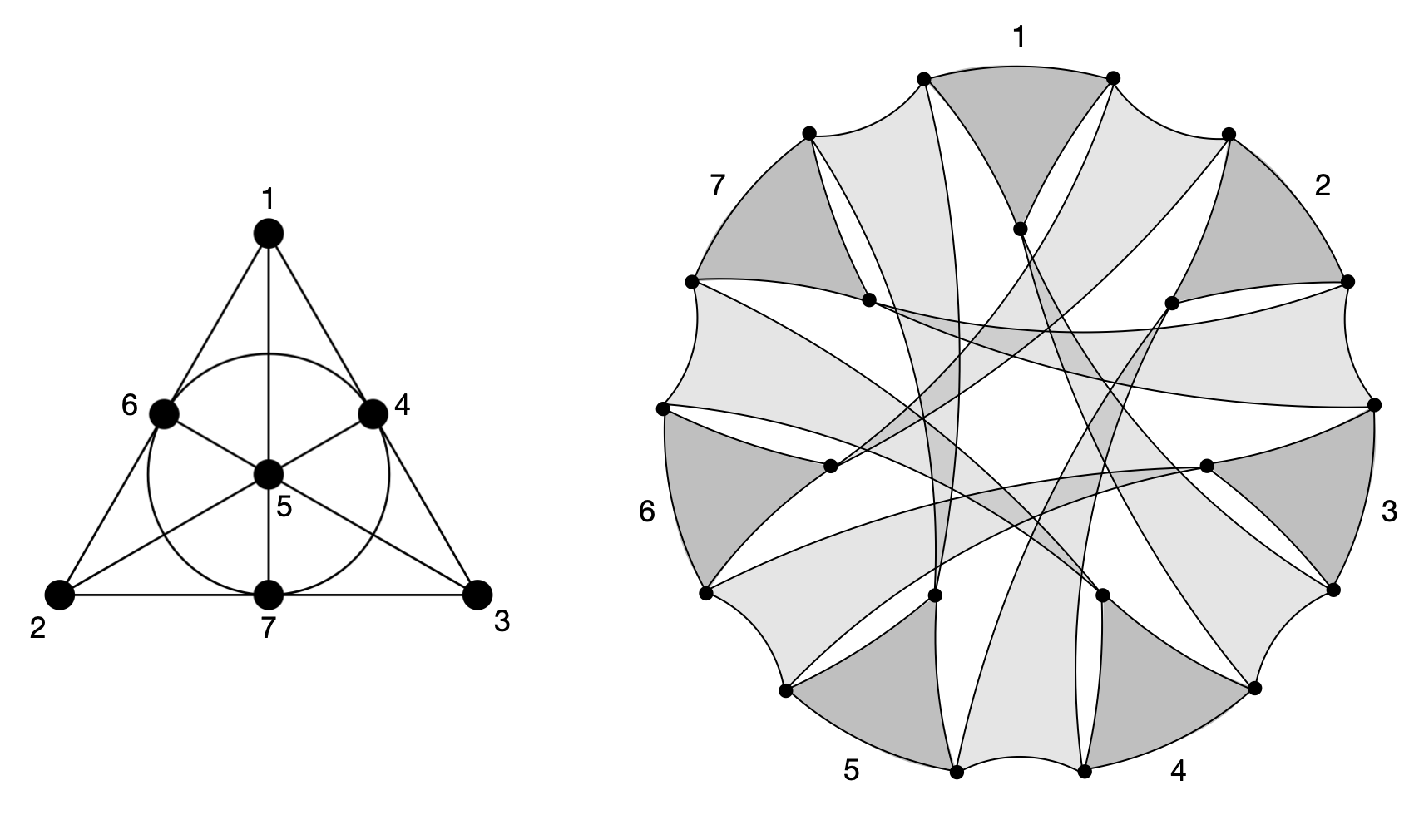}
    \caption{Example for Theorem \ref{t:small_Steiner}, case $t = 2$, $k = 3$: the Steiner system $S(2,3,7)$ is given by the Fano plane $F$ (left), the graph $S_{2}(3) = \Gamma(K_3,F)$ is on the right.}
    \label{f:ex_Steiner}
\end{figure}

In the following, we apply Theorem \ref{t:general_method} to present a recursive construction of a family $A_t(q)$ of graphs that employs affine geometries. 
It is well known that, given a prime power $q$, there exists a Steiner system $S(2, q, q^m)$, for any $m \ge 2$, that can be extracted from the $m$-dimensional affine geometry $AG_m(q)$ over $\mathbb{F}_q^m$ by taking the vector space $X = \mathbb{F}_q^m$ and all $1$-dimensional affine subspaces (i.e., lines), that represent the blocks, see \cite{Sti04}. Such a system has the property of being \emph{resolvable}, meaning that the set of blocks can be partitioned into $r = \frac{q^m - 1}{q-1}$ parallel classes $\cB_1, \cB_2, \ldots, \cB_r$, where each $\cB_i$ contains $q^{m-1}$ pairwise disjoint blocks whose union is $X$, $1 \le i \le r$.

\btm \label{t:affine_geom}
Let $t$ and $q$ be positive integers, where $q$ is a prime power.  Let $n_t = q^{2^t-1}$, for $t \ge 1$. Let $H_t(q)$ be a subhypergraph of the affine geometry $AG_{2^t}(q)$ that is induced by the hyperedges of $n_t$ arbitrarily chosen parallel classes. Define 
\[A_1(q) = K_q \mbox{ and } A_{i+1}(q) = \Gamma(A_i(q), H_i(q)), \;i \geq 1.\]
Then the graph $A_t(q)$ has order $n_t$ and constant link $tK_{q-1}$, $t \geq 1$. Moreover, the graph $A_t(3)$ is $\{C_4,C_5\}$-free and the graph $A_t(4)$ is $C_5$-free, for every $t \ge 1$.
\etm

\bpf
Since the affine geometry $AG_{2^t}(q)$ has $v_t = q^{2^t}$ points, blocks of size $q$, and $\frac{v_t-1}{q-1} = \frac{q^{2^t}-1}{q-1} \ge q^{2^t-1} = n_t$ parallel classes, the hypergraph $H_t(q)$ is well defined and it is $n_t$-regular and $q$-uniform. By Theorem \ref{t:general_method}, the graph $A_t(q)$ has constant link $tK_{q-1}$, for every $t \geq 1$. Moreover, the same theorem implies that $A_t(q)$ has order 
\[q \prod_{i=1}^{t-1}v_i = q \prod_{i=1}^{t-1} q^{2^i} = q^{\sum_{i= 0}^{t-1} 2^i} = q^{2^t-1} = n_t.\]
As in $H_t(q)$ every pair of vertices is contained in exactly one hyperedge, $g(H_t(q)) \geq 3$ holds. Moreover, $A_1(q)$ has circumference $q$. Hence, we can apply item (ii) of Theorem \ref{t:general_method} for the cases where $q \leq 5$. It follows that $A_t(3)$ is $\{C_4,C_5\}$-free, and $A_t(4)$ is $C_5$-free, for every $t \ge 1$.
\epf

\begin{figure}[ht]
    \centering
    \includegraphics[width=0.7\linewidth]{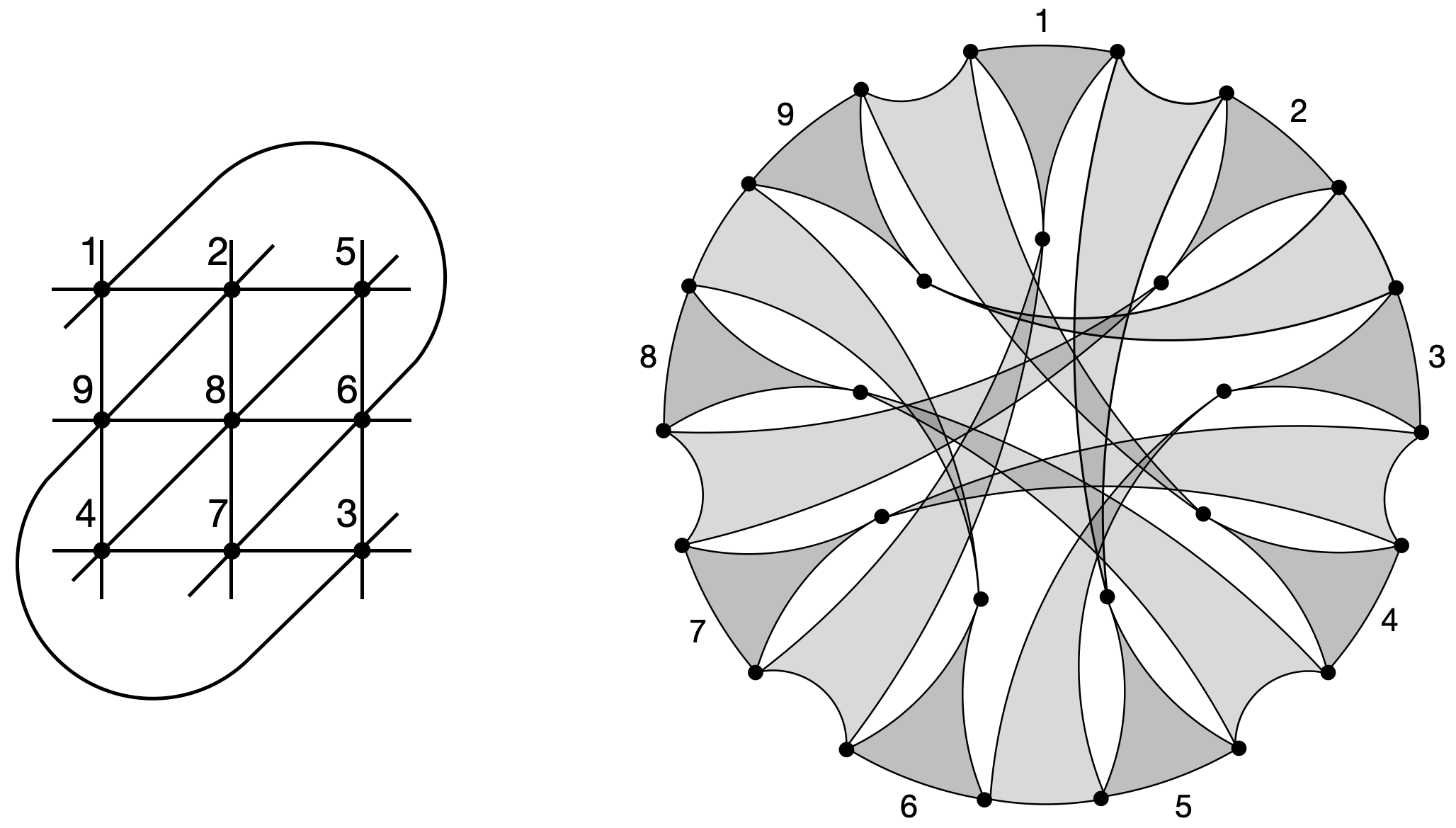}
    \caption{Example for Theorem \ref{t:affine_geom}, case $t = 2$, $q = 3$: three parallel classes of the affine geometry $AG_2(3)$ are depicted left, that are used to construct the graph $A_{2}(3) = \Gamma(K_3,H_1(3))$ on the right.}
    \label{f:ex_affine_geom}
\end{figure}

\section{Connected graphs}\label{s:connected_graphs}

In this section, we deal with the problem of providing connected graphs of infinitely many orders that attain equality in Theorem \ref{t:tur_ind_Ck_star}. As already mentioned before,  we may assume that the hypergraph given by Theorem \ref{t:existence_Furedi_etal} is connected (just by taking 
one connected component in a disconnected case).

Going in \cite{ElLin14} into the proof of the second existence theorem we made reference to, Theorem \ref{t:d-r}, we can fix the parameters $k$ and $t$ and let $g$ grow to convince oneself that there are infinitely many connected hypergraphs satisfying these conditions. Hence, this provides us already with an infinite family of graphs for Theorem \ref{t:existence_const_link-gr}, by means of which we proved the sharpness of Theorem \ref{t:tur_ind_Ck_star}. However, Theorem \ref{t:d-r} gives an upper bound on the order of these graphs that grows exponentially in $g$. It similarly occurs with Theorem \ref{t:sep_perm_constr}: we can construct connected graphs with the required properties and as large as desired just by taking a sufficiently large $q$ (the number of iterations in Construction \ref{const:k-ary-tree-like}). Again, the number of vertices grows exponentially with $q$. We will see here that the order can be controlled to a linear grow---in fact, forming an arithmetic progression---by means of very simple operations on the desired number of copies of just one hypergraph.

There are many operations in graph and hypergraph theory called switching, mainly with the meaning of replacing two edges by another pair of edges. The most familiar one is the following: given a graph $G$ and four vertices $a,b,x,y$ where $ab$ and $xy$ are edges in $G$ and $xb$ and $ya$ are non-edges, removing $ab$ and $xy$ and inserting $xb$ and $ya$ is called \emph{switching}. This is a basic operation that has many applications in degree sequences of graphs as it preserves degrees. For recent papers about switching we refer to \cite{BEFHRST13, Fer22} and the references therein. A generalization of switching to $3$-uniform hypergraphs is given in \cite{HMBBM}.

Borrowing from the common practice of using switching to enlarge triangle-free regular graphs (and large girth in general), we introduce a general switching operation in hypergraphs in order to create larger connected $t$-regular $k$-uniform hypergraphs with large girth from smaller ones. Let $H$ be a $k$-uniform hypergraph. Let $e, f \in E(H)$ be two edges, and let $e = e_1 \cup e_2$ and $f = f_1 \cup f_2$ be partitions of them such that $0 < |e_1| = |f_1| < k$. Let $e^* = e_1 \cup f_2$ and $f^* = e_2 \cup f_1$. We transform the hypergraph $H \mapsto (H - \{e, f\}) + \{e^*, f^*\}$ as the one obtained from $H$ by switching $e$ and $f$ with $e^*$ and $f^*$.

\begin{figure}[ht]
    \centering
    \includegraphics[width=0.95\linewidth]{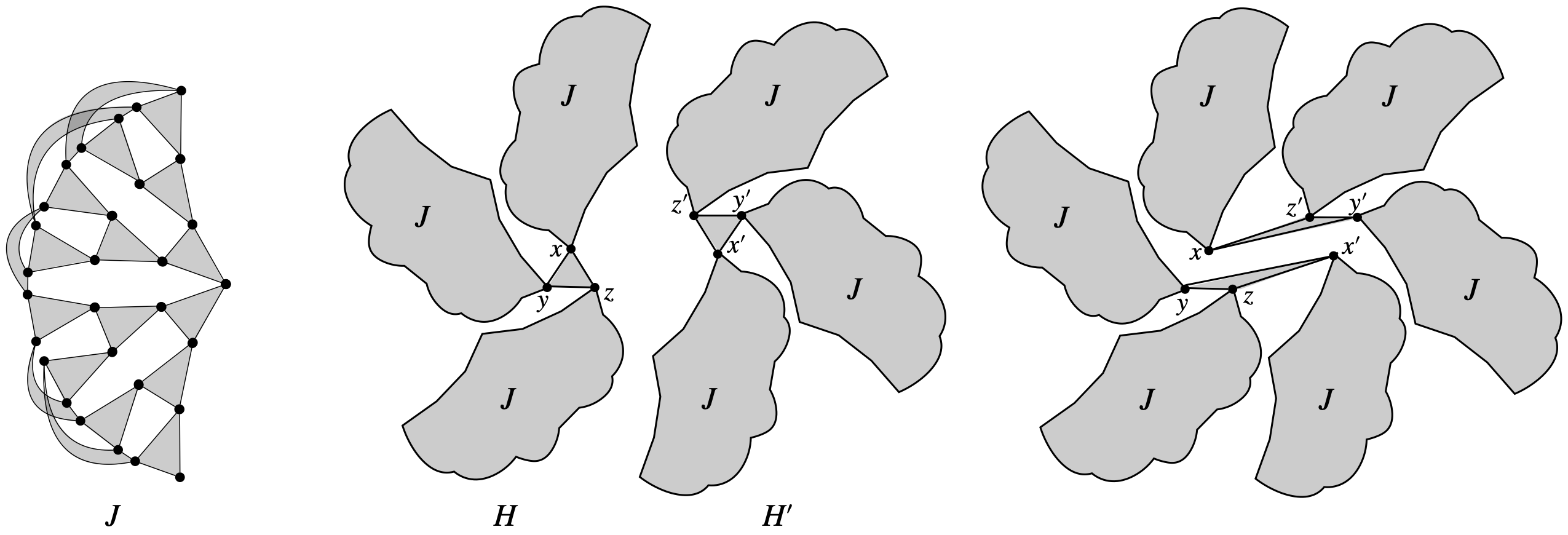}
    \caption{Hypergraph $H$ is built from three copies of the $2$-regular $3$-uniform hypergraph $J$, together with one hyperedge containing the three degree-$1$ vertices $x, y, z$. Take a copy $H'$ of $H$, with vertices $x', y', z'$ being the copies of $x, y ,z$. The hypergraph obtained from $H \cup H'$ by switching $\{x,y,z\}, \{x', y', z'\}$ with $\{x,y',z'\}, \{x', y, z\}$ yields again two separate copies of $H$, hence a disconnected hypergraph.}
    \label{f:disconn_switching}
\end{figure}

We observe at this point that applying the switching operation on two connected $k$-uniform $t$-regular hypergraphs with respect to any choice of hyperedges and vertices does not necessarily yield a connected hypergraph. An example is exhibited in Figure \ref{f:disconn_switching}.

We will show that applying the switching operation in an appropriate manner will yield larger connected hypergraphs from smaller ones. Before proceeding, we need to prove the following basic property of connected hypergraphs.

\begin{lemma}\label{l:non-cut-edges}
   Let $H$ be a connected hypergraph without trivial edges and with $\delta(H) \geq 2$. Then $H$ contains an edge whose removal results in a connected hypergraph.
\end{lemma}

\bpf
The removal of a cut-edge from $H$ leaves a vertex partition into classes that induce non-extendable connected subhypergraphs (components). Since $\delta(H) \geq 2$, all those parts contain edges of $H$. Among all edges of $H$, choose an edge $e_0$ whose removal yields a smallest part. Say $H-e_0$ has components with vertex sets $X_1, X_2, \ldots, X_q$, and let $X_1$ have the smallest size. Let $e_1$ be an edge inside $X_1$. In $H-e_1$, due to the presence of $e_0$, the set $X_2 \cup \cdots \cup X_q$ belongs to one big component together with some further vertices of $X_1$, at least those of $e_0$ belonging to $X_1$. So, if $e_1$ was a cut-edge, it would yield a component of an order smaller than $|X_1|$, contradicting the choice of $e_0$.
\epf

\begin{prop}[Switching]\label{p:switching} Let $H_1, H_2$ be two connected and disjoint $k$-uniform $t$-regular hypergraphs of respective girths $g_1, g_2 \geq 3$, where $t, k \geq 2$. 
Then there are edges $e \in E(H_1)$ and $f \in E(H_2)$, such that, if $e = e_1 \cup e_2$ and $f = f_1 \cup f_2$ are partitions with $0 < |e_1| = |f_1| < k$, and $e^* = e_1 \cup f_2$ and $f^* = e_2 \cup f_1$, then the hypergraph obtained from $H_1 \cup H_2$ by switching $e$ and $f$ by $e^*$ and $f^*$ is a connected $k$-uniform $t$-regular hypergraph of girth at least $g = \min \{g_1, g_2\}$. 
\end{prop}

\bpf
As given by Lemma \ref{l:non-cut-edges}, we can choose edges $e \in E(H_1)$ and $f \in E(H_2)$ such that both $H_1 - e$ and $H_2 - f$ are connected. Let $e = e_1 \cup e_2$ and $f = f_1 \cup f_2$ be partitions such that $0 < |e_1| = |f_1| < k$. Now, inserting the edges $e^* = e_1 \cup f_2, f^* = e_2 \cup f_1$ to the union $(H_1 - e) \cup (H_2 - f)$ yields a connected hypergraph $H$, which is clearly $k$-uniform, and is the one obtained from $H_1 \cup H_2$ by switching $e$ and $f$ by $e^*$ and $f^*$.

Consider now any vertex $v \in V(H)$. Observe that $W:= e \cup f = e^* \cup f^*$. If $v$ does not belong to $W$, then its degree is the same as before. If $v \in W$, then it belongs to one edge less from the original hypergraph ($H_1$ or $H_2$) and to one new edge in the new hypergraph $H$. Hence, the degree remains the same, and so $H$ is $t$-regular.

Consider now a cycle $C$ of $H$. Assume first that $C$ has none of the edges $e^*$, $f^*$. Then it is contained either in $H_1-e$ or in $H_2 - f$, and so it has length at least $g = \min \{g_1,g_2\}$. Now assume that $C$ has some edge of $e^*$ or $f^*$. Then it necessarily has to contain both, as otherwise one cannot close the cycle. Since the edges $e^*$, $f^*$ are disjoint, the cycle $C$ has to go, inside $H_1-e$ through a path $P_1$ joining some vertex in $e_1$ with some vertex in $e_2$. The path $P_1$ together with $e$ makes a cycle in $H_1$, which has length at least $g_1$, implying that $P_1$ has length at least $g_1-1$. Similarly occurs with the part of the cycle that is inside $H_2-f$: it is a path $P_2$ of length at least $g_2 - 1$. Hence, $C$ has length at least $g_1 + g_2 > g$. It follows that $H$ has girth at least $g = \min \{g_1, g_2\}$.
\epf

\btm \label{t:connected}
Let $t, k, g$ be integers such that $t, k \geq 2$, and $g \geq \max\{4, k+1\}$, and let $F$ be a graph of order $k$ and constant link $L$ and with circumference $c$. Then there is an integer $n$ such that, for any $q \geq 1$, there is a connected graph of constant link $tF$ and order $qn$ without cycles of length $\ell$ for any $c < \ell < g$. 
\etm

\bpf
We will show first, by induction on $q \geq 1$, that there is a connected $t$-regular $k$-uniform hypergraph $H_q$ of girth at least $g$ and order $qn$, for some integer $n$.

To this aim, let $G$ be a graph with constant link $tK_{k-1}$ free of cycles of length $\ell$ for any $k < \ell < g$, whose existence is guaranteed by Theorem \ref{t:existence_const_link-gr}. By Proposition \ref{p:coloring} and Theorem \ref{t:families_equiv}, we can transform $G$ into a $t$-regular $k$-uniform hypergraph $\varphi^{-1}(G) = H_1$ of girth at least $g$ by replacing every copy of $K_k$ with a hyperedge. Setting $n = n(G)$, we are done with the induction start $q = 1$.

Suppose now there is a connected $t$-regular $k$-uniform hypergraph $H_q$ of girth at least $g$ and order $qn$. We can now apply the switching operation of Proposition \ref{p:switching} on $H_1$ and $H_q$ to obtain a connected $t$-regular $k$-uniform hypergraph $H_{q+1}$ of girth at least $g$ and order $(q+1)n$, as desired.

To finish the proof, we can use back on Theorem \ref{t:families_equiv} to build a graph $\varphi(H_q)$ of order $qn$ by replacing every hyperedge of $H_q$ by a copy of $K_k$ and then apply Proposition \ref{p:coloring} to replace every $K_k$ by a copy of $F$. The graph obtained in this way $G'$ has no cycles of length $\ell$ for any $c < \ell < g$ by Proposition \ref{p:coloring}.
\epf

Proposition \ref{p:switching} allows a further alternative proof, via switching, of the main result of this paper.

\btm \label{t:connected_linear_sequence}
Let $t,k,g$ be integers such that $t, k \geq 2$ and $g \geq \max\{4, k+1\}$, and let $F$ be a connected graph of order $k$ and  constant link $L$ and with circumference $c$. Then there is an integer $n$ such that, for any $q \geq 1$, there is a connected graph of constant link $tL$ and order $qn$ without cycles of length $\ell$ for any $c < \ell < g$. 
\etm

\bpf
Let $G$ be a connected graph with constant link $tK_{k-1}$ free of cycles of length $\ell$ for any $k < \ell < g$, whose existence is guaranteed by Theorems \ref{t:existence_const_link-gr} and \ref{t:sep_perm_constr}. Let $n$ be the order of $G$. By Proposition \ref{p:coloring} and Theorem \ref{t:families_equiv}, we can transform $G$ into a connected $t$-regular $k$-uniform hypergraph $H = \varphi^{-1}(G)$ of girth $g$ by replacing every copy of $K_k$ with a hyperedge. This settles the case where $q = 1$. For $q \geq 2$, take $q$ copies of $H$, say $H_1, H_2, \ldots, H_q$, and perform the switching operation given in Proposition \ref{p:switching} to obtain a connected $t$-regular $k$-uniform hypergraph $H'$ of girth at least $g$ and order $qn$. Now we can use back again Theorem \ref{t:families_equiv} and Proposition \ref{p:coloring} to build a graph $G'$ of order $qn$ by replacing every hyperedge by a copy of $F$. By Proposition \ref{p:coloring}, $G'$ is connected and has no cycles of length $\ell$ for any $c < \ell < g$.
\epf

As already mentioned in the Introduction, it is important to remark that Theorem \ref{t:connected_linear_sequence} does not work for $t = 1$ in general, since there are examples of graphs that appear as link graphs of only finitely many connected graphs, as is $K_k$ and the Petersen graph \cite{Hal80}.

A very important consequence of Theorem \ref{t:connected_linear_sequence} is that the induced Turán numbers attaining the upper bound in Theorem \ref{t:tur_ind_Ck_star} occur for all sufficiently large $n \equiv 0 \pmod{n(G)}$ via connected graphs, just by taking as many copies of $G$ as needed, where $G$ is a single constructed example of constant link $tK_{k-1}$, and then applying the switching operation.

\begin{cor}
    The bound in Theorem \ref{t:tur_ind_Ck_star} is attained by an infinite collection of connected graphs, whose orders form an arithmetic progression.
\end{cor}

Finally, by means of Theorem \ref{t:families_equiv}, we can state the following corollary to Theorem \ref{t:connected_linear_sequence}.
\begin{cor}
For any integers $t,k,g$ such that $t, k \geq 2$, and $g \geq \max\{4, k+1\}$, the families $\mathcal{H}_{t,k,g}$, $\mathcal{L}_{t,k,g}$ and $\mathcal{B}_{t,k,g}$ contain an infinite collection of connected members, whose orders form an arithmetic progression. 
\end{cor}

\paragraph{Acknowledgements.}

We gratefully thank Raphael Yuster for his kindly providing us with a preliminary version of Lemma \ref{l:sep_perm}. Research of the second author was supported by grant PAPIIT IG100822. Research of the third author was supported in part by the National Research, Development and Innovation Office, NKFIH Grant FK 132060.


\begin{thebibliography}{9}

\bibitem{AdBryBu07}
P. Adams,  D. Bryant,  M. Buchanan,
A survey on the existence of $G$-Designs,
J. Combin. Designs 16 (2008), 373--410.

\bibitem{AJRS22}
G. Araujo-Pardo, R. Jajcay, A. Ramos-Rivera, T.  Sz\H{o}nyi,
On a relation between bipartite biregular cages, block designs and generalized polygons, J. Combin. Des. 30(7) (2022), 479--496.

\bibitem{BEFHRST13}
S. Behrens, C. Erbes, M. Ferrara, S. G. Hartke, B. Reiniger, H. Spinoza, C. Tomlinson,
New results on degree sequences of uniform hypergraphs,
Electron. J. Combin. 20 (2013), no. 4, Paper 14, 18 pp.


\bibitem{Bla23}
S. R. Blackburn,
Permutations that separate close elements,
J. Combin. Th. (A) 196 (2023), 105734.


\bibitem{BlaEt23}
S. R. Blackburn, T. Etzion
Permutations that separate close elements, and rectangle packings in the torus,
arXiv:2306.03685 (2023).



\bibitem{BlaHM80}
A. Blass, F. Harary, Z. Miller,
Which trees are link graphs?,
J. Combin. Th. (B) 29 (1980) 277--292.

\bibitem{Brou23}
A. E. Brouwer, 
Some locally Kneser graphs, 
arXiv:2312.02964 (2023).
%https://www.win.tue.nl/~aeb/preprints/loc-kneser.pdf

\bibitem{Brou95}
A. E. Brouwer, 
Variations on a theme by Weetman,
14th British Combinatorial Conference (Keele, 1993),
Discrete Math. 138(1--3) (1995), 137--145.
%https://core.ac.uk/reader/82665920

\bibitem{BroCon73}
M. Brown, R. Connelly,
On graphs with a constant link,
pp. 19--51 in:
New Directions in the Theory of Graphs, Academic Press, New York, 1973.

\bibitem{BroCon75}
M. Brown, R. Connelly,
On graphs with a constant link, II,
Discrete Math. 11 (1975) 199--232
%https://www.sciencedirect.com/science/article/pii/0012365X75900370


\bibitem{Bul73}
V. K. Bulitko,
On graphs with given vertex-neighbourhoods,
Trudy Mat. Inst. Im. Steklova 133 (1973) 78–94.
(Russian)

\bibitem{Chv77}
V. Chvátal, 
Tree-complete graph Ramsey numbers,
J. Graph Theory 1 (1977), no.1, 93.


\bibitem{CD-hand}
C. J. Colbourn, J. H. Dinitz (Eds.), Handbook of Combinatorial Designs, 2nd ed. Discrete Mathematics and Its Applications, Chapman and
Hall/CRC, Boca Raton, FL, USA, 2006.


\bibitem{ElLin14}
D. Ellis, N. Linial,
On regular hypergraphs of high girth.
Electron. J. Combin. 21(1) (2014), \#P1.54

\bibitem{ErSa63}
P. Erd\H{o}s, H. Sachs,
Regul\"are Graphen gegebener Taillenweite mit minimaler Knotenzahl (German),
Wiss. Z. Martin-Luther-Univ. Halle-Wittenberg Math.-Natur. Reihe 12(1963), 251--257.

\bibitem{EGM19} B. Ergemlidze , E. Gy\H{o}ri, A. B. Methuku,  
Turán Number of an Induced Complete Bipartite Graph Plus an Odd Cycle, Combinatorics, Probability and Computing, 28(2) (2019), 241--252.

\bibitem{ErTu23}
G. Erskine, J. Tuite,
Small graphs and hypergraphs of given degree and girth,
Electron. J. Combin. 30(1) (2023), Paper No. 1.57, 8 pp.


\bibitem{FLSUW95}
Z. F\H{u}redi, F. Lazebnik, Á. Seress,  V. A. Ustimenko, A. J. Woldar, 
Graphs of prescribed girth and bi-degree, J. Combin. Theory Ser. B, 64(2) (1995), 228--239.

\bibitem{FaScSh80}
R. J. Faudree, R. H. Schelp, J. Sheehan, Ramsey numbers for matchings,
Discrete Math. 32 (1980), 105--123.


\bibitem{Fer22}
R. Fernandes,
On the switch-length of two connected graphs with the same degree sequence
Australas. J. Combin. 83 (2022), 87--100.


\bibitem{Fro89}
D. Fron\v{c}ek,
Locally linear graphs,
Math. Slovaca 39 (1989), no.1, 3--6.
%https://dml.cz/bitstream/handle/10338.dmlcz/136481/MathSlov_39-1989-1_1.pdf 

\bibitem{Hal80}
J. I. Hall,
Locally Petersen Graphs,
J. Graph Th. 4 (1980), 173--187.


\bibitem{Hal85}
J. I. Hall,
Graphs with constant link and small degree or order,
J. Graph Th. 8 (1985) 419--444.


%\bibitem{Han61}
%H. Hanani, 
%The existence and construction of balanced incomplete block designs, 
%Ann. Math. Stat. 32 (1961), 361--386.


%\bibitem{Haw23}
%E. Hawboldt, 
%A machine learning approach to constructing Ramsey graphs leads to the Trahtenbrot-Zykov problem,
%Electronic Theses and Dissertations (2023), Paper 4125.


%\bibitem{HoWo16}
%C. Hoppen, N. Wormald,
%Properties of regular graphs with large girth via local algorithms,
%J. Combin. Theory Ser. B 121 (2016), 367--397.
%https://www.sciencedirect.com/science/article/pii/S0095895616300569

\bibitem{Ill21-1}
F. Illingworth,
A note on induced Turán numbers,
arXiv:2105.12503.

\bibitem{Ill21-2}
F. Illingworth,
Graphs with no induced $K_{2,t}$,
Electron. J. Combin. 28 (1) (2021), \# P1.19.


\bibitem{Kee14}
P. Keevash,
The existence of designs,
arXiv:1401.3665.

%\bibitem{Kirk}
%T. P. Kirkman, 
%On a problem in combinations, 
%Cambridge and Dublin Math. J. 2 (1847), 191--204

\bibitem{HMBBM}
A. Hubai, T.R. Mezei, F. Béres, A. Benczúr, I. Miklós, 
Constructing and sampling partite, $3$-uniform hypergraphs with given degree sequence,
PLoS ONE 19(5) (2024), e0303155.


\bibitem{LaPiVi19}
F. Larrión, M. A. Pizaña, R. Villarroel-Flores, 
On the clique behavior of graphs with small constant link,
Ars Combin. 142 (2019), 27--53.
%http://xamanek.izt.uam.mx/map/papers/locally6.pdf 


\bibitem{LaMu13}
R. C. Laskar, H. M. Mulder, Henry Martyn,
Path-neighborhood graphs,
Discuss. Math. Graph Theory 33 (4) (2013), 731--745.
%https://www.dmgt.uz.zgora.pl/publish/pdf.php?doi=1700 

\bibitem{LaSc16}
J. Lauri, R. Scapellato,
Topics in graph automorphisms and reconstruction,
London Math. Soc. Lecture Note Ser. 432,
Cambridge University Press, Cambridge (2016) xiv+192 pp.


\bibitem{LTTZ18}
P.-S. Loh, M. Tait, C. Timmons, R. M. Zhou,
Induced Turán Numbers, 
Combinatorics, Probability and Computing, 27(2) (2018), 274--288.

\bibitem{NTT18}
V. Nikiforov, V., M. Tait, C. Timmons,
Degenerate Turán Problems for Hereditary Properties,
Electron. J. Combin. 25(4)(2018) \# P4.39.

\bibitem{Ofi19}
D. Ofir, 
Graphs with large girth and free groups,
arXiv:1907.06936 (2019).
%https://arxiv.org/pdf/1907.06936.pdf


\bibitem{Rad_survey}
S. P. Radziszowski, 
Small Ramsey numbers,
Electron. J. Combin. 1 (1994), Dynamic Survey 1, 30 pp.

%\bibitem{Sko58}
 %T. Skolem, Some remarks on the triple systems of Steiner, Math. Scand. 6 (1958), 273--280.


\bibitem{Sti04}
D. R. Stinson,
Combinatorial designs. Constructions and analysis.
Springer-Verlag, New York (2004), xvi+300 pp.


\bibitem{Toma89}
J. Tomanová,
A note on link graphs,
Math. Slovaca 39(3) (1989), 225--231.
%https://dml.cz/bitstream/handle/10338.dmlcz/129354/MathSlov_39-1989-3_1.pdf
%   clique-star in link.

\bibitem{Viz65}
V. G. Vizing, 
Critical graphs with given chromatic class, 
Metody Diskret. Analiz. 5 (1965), 9--17.


\bibitem{Vog86}
W. Vogler, 
Representing groups by graphs with constant link and hypergraphs,
J. Graph Th. 10 (1986), 461--475.
%https://www.semanticscholar.org/paper/Representing-groups-by-graphs-with-constant-link-Vogler/4fda0555855a26ef64c53968ac63a3e7deb83dcb 

\bibitem{Weet94-1}
G. M. Weetman, 
A construction of locally homogeneous graphs, 
J. London Math. Soc. (2) 50 (1) (1994), 68--86. 

\bibitem{Weet94-2}
G. M. Weetman, 
Diameter bounds for graph extensions, 
J. London Math. Soc. (2) 50 (2) (1994), 209--221.

\bibitem{BookWest}
D. B. West,
Introduction to graph theory,
Prentice Hall, Inc., Upper Saddle River, NJ (1996), xvi+512 pp.

\bibitem{Wil75}
R. M. Wilson, 
Decompositions of complete graphs into subgraphs isomorphic to a given graph, 
Proc. British Combinatorial Conf. (1975), 647--659.

\bibitem{YuWu21}
H. Yu, B. Wu, 
Graphs in which $G-N[v]$ is a cycle for each vertex $v$,
Discrete Math. 344(2021), no.9, Paper No. 112519, 7 pp.

\bibitem{Zel86}
B. Zelinka,
Edge neighbourhood graphs,
Czechoslovak Math. J. 36 (111) (1986), no.1, 44--47.
%https://dml.cz/bitstream/handle/10338.dmlcz/102064/CzechMathJ_36-1986-1_7.pdf

\bibitem{Zyk63}
A. A. Zykov,
Problem 30, pp. 164--165. in:
Theory of graphs and its applications,
Proc. Symp. Smolenice 1963 (M. Fiedler, ed.),
Academic Press, Prague 1964.


\end{thebibliography}
\end{document}